# Theoretical Analysis on the Ability of Particle Swarm Optimization to Escape from Local Optimum

Haoxin Wang[a,b], Libao Shi[b*]

[a] Department of Electrical Engineering, Tsinghua University, Beijing, 100084, P. R. China

[b] Institute for Ocean Engineering, Tsinghua Shenzhen International Graduate School, Tsinghua University, 518055, P. R. China

ABSTRACT

As one of the most prominent swarm intelligence algorithms, particle swarm optimization (PSO) has been extensively applied to solve global optimization problems. The theoretical analysis on the ability of PSO to escape from local optimum (LO) is deeply intertwined with its global optimization performance, making it highly significant to investigate the underlying mechanisms that enable PSO to escape from LO and achieve global optimization. In this study, a rigorous mathematical analysis is presented to explore the ability and mechanisms of PSO to escape from LO. By modelling each search agent as a Markov chain, the probability of PSO escaping from LO within a finite number of iterations is calculated, and the behavior of agents prior to escaping from LO is analyzed. Through this theoretical analysis, two necessary and sufficient conditions are proposed to ensure that PSO escapes from LO in finite iterations with a probability 1: (a) the personal best (pbest) and global best (gbest) positions stagnate around different local optima, and (b) the inertial weight in PSO is equal to one. Upon these conditions, two potential behaviors of agents, namely inertial motion and oscillation, are proposed, and a state transition chain of the agent toward the global optimum is constructed. Additionally, the theoretical analysis reveals that increasing the distance between pbest and gbest, as well as the values of the learning factors in PSO, enhances the speed at which PSO escapes from LO. Finally, the conclusions drawn from the theoretical analysis are corroborated through numerical experiments.

## 1. INTRODUCTION

In recent years, a variety of swarm intelligence (SI) algorithms have emerged as effective tools for addressing increasingly complex optimization problems in real-world applications [1-2]. Unlike traditional gradient-based optimization methods, such as the primal-dual interior-point method [3] and the Alternating Direction Method of Multipliers (ADMM) [4], SI algorithms are characterized by two key features: (1) the existence of multiple search agents, where the position update of each agent depends on the positions of certain other agents, and (2) the integration of randomness into the position update mechanism, ensuring that new positions of agents follow a specific probability distribution. These features enable SI algorithms to explore the entire search space more effectively, reducing the likelihood of becoming trapped in local optima, which is a common drawback of traditional optimization methods [5]. Consequently, SI algorithms are frequently employed for solving global optimization problems. Among these, Particle Swarm Optimization (PSO) stands out as one of the most widely recognized and applied algorithms [6]. Owing to its straightforward and transparent position update mechanism, as well as its robust optimization performance, PSO has been successfully utilized in diverse real-world applications, including economic dispatch [7-8], path planning [9-10], neural network training [11-12], and feature selection [13]. To further enhance PSO's performance, extensive studies have focused on improving various aspects of the algorithm, such as parameter selection [14-16], strategies for selecting particle learning agents [17-19], velocity or position update mechanism [20-23], and the incorporation of subpopulations [24].

Alongside the extensive research on applications and improvements of PSO, significant efforts have also been directed toward its theoretical analysis. As previously highlighted, a key advantage of PSO and other SI algorithms lies in their robust global search capability. Consequently, a fundamental aspect of PSO's theoretical study is its global convergence. As a stochastic search algorithm, an essential metric for evaluating PSO's global convergence is the probability of finding the Global Optimum (GO) within a finite number of iterations. If this probability equals one, the algorithm is said to find the GO with probability 1, indicating strong global convergence. Early work in [25] established a sufficient condition to ensure that a stochastic search algorithm finds the GO with probability 1. This foundation work paved the way for more specific analyses of PSO. For instance, the PSO's velocity and position update mechanism was analyzed [26], demonstrating through case studies that PSO may not find the GO with probability 1. Subsequently, it was proved that PSO is global convergent only under certain conditions using probability theory [27]. Similarly, the concept and properties of martingales were employed in stochastic processes to show that PSO finds the GO with a probability less than 1 [28]. More recently, researchers have adopted alternative approaches by formulating the particle dynamics as a system of stochastic differential equations in continuous time. Techniques such as the Laplace principle were then applied to analyze the equilibrium points of PSO and their proximity to the GO [29-30]. However, existing research based on stochastic processes often provides only a macroscopic and somewhat generalized understanding, while the stochastic differential equation-based approach remains highly abstract. As a result, these findings offer limited insights into the actual operational mechanisms of PSO for finding the GO. Key questions remain unresolved: What behavior patterns do agents exhibit prior to reaching the GO? How do the initial positions of agents and the parameter selection mechanisms in PSO influence the probability of finding the GO within a finite number of iterations? These specific issues necessitate further in-depth theoretical analysis.

Since the original PSO may lack global convergence, another significant research direction has focused on designing PSO variants with guaranteed global convergence. These studies are often motivated by the following observation: Given the positions and velocities of agents at iteration $t$, the velocities at iteration $t+1$ can be expressed as the sum of two uniformly distributed random variables (see PSO's velocity update equation in Section 2.1). This implies that the position distribution at $t+1$ is bounded, potentially failing to cover the entire feasible region. Consequently, there exist iterations where PSO lacks global search capability, leading to the absence of a theoretical guarantee for PSO's global convergence. To address this limitation, researchers have modified PSO's velocity and position update mechanisms to follow an unbounded distribution or at least a distribution that covers the entire feasible space. This ensures that, in every iteration, each agent always has a positive probability of finding the GO, thereby allowing the algorithm to find the GO with probability 1. For instance, studies in [21,31-34] introduced normally distributed random variables into the position update formula. Quantum PSO [35-37] employed the logarithm transformation of uniformly distributed random variables to construct positions. Other approaches have utilized random variables following Lévy [18,38-39] or Cauchy distributions [40] to define agent velocities and positions, while some researchers have incorporated chaotic numbers generated by discrete dynamical systems into the velocity update expression [41-42]. Although these modifications ensure the global convergence of the algorithm, they do not fundamentally resolve the underlying issue with PSO's global convergence. Specifically, while the agent's distribution at iteration $t+1$ may not cover the entire feasible region, the distributions at iterations $t+2, t+3$, and beyond remain uncertain, leaving the global searching capability of the agent in subsequent iterations unclear. In summary, the global optimization capability and mechanisms of PSO still require further in-depth investigation.

To provide a clearer understanding of the global optimization mechanism of PSO, this study focuses on a specific and fundamental issue: *the ability of PSO to escape from Local Optimum (LO)*. We believe that investigating PSO's capability to escape from LO contributes to a more detailed and fundamental characterization of its global convergence. This perspective is supported by the following two key considerations:

(1) The ability to escape from LO is a necessary condition for an algorithm to achieve global convergence. Extensive research has shown that SI algorithms solve global optimization problems through two primary steps [52-56]: first, agents explore the entire feasible region with relatively large step sizes to identify a rough region containing the GO, a process referred to as *exploration*; second, agents exploit the identified region with smaller step sizes to refine the solution and obtain the precise optimal value, a process known as *exploitation*. However, for complex global optimization problems characterized by multiple local optima or highly deceptive local optima (e.g., a LO with a wide neighborhood or a significantly small objective function value, as discussed in [57]), the region identified during the exploration step is more likely to be the neighborhood of a LO rather than the GO. Consequently, after completing these two steps, the algorithm may fail to achieve global optimization and instead become trapped in a LO. Under such circumstances, only if the algorithm possesses the ability to escape from LO can it continue the global optimization process. Therefore, the ability to escape from LO is a necessary condition for ensuring the global convergence of an algorithm.

(2) An algorithm capable of escaping LO inherently exhibits global convergence. Compared to other population states, when trapped in a LO, the population becomes highly concentrated, the search range of agents is significantly reduced, and the positions of guiding agents (i.e., those with higher fitness, such as the pbest and gbest in PSO, as described in Section 2.1) become more difficult to update. If the population can escape from the LO under such conditions and continue the global optimization process, it can be inferred that the algorithm possesses robust global optimization capabilities, regardless of the population's state. Therefore, the ability to escape from LO serves as a sufficient condition for the global convergence of an algorithm.

This study presents a theoretical analysis framework to investigate the ability and mechanism of PSO to escape from LO. By jointly modelling the position and velocity of each search agent as a continuous state Markov chain, we rigorously calculate the probability of PSO escaping from LO within a finite number of iterations and analyze the behavior of agents prior to escaping from LO. Theoretical analysis reveals that two necessary and sufficient conditions must be satisfied to ensure that PSO escapes from LO in finite iterations with a probability of 1: (1) the *pbest* and *gbest* (see Section 2.1) positions stagnate around different local optima, and (2) the inertia weight $\omega$ is equal to one. Under these conditions, two possible behaviors of agents, namely inertial motion and oscillation, are identified, and a state transition chain leading toward the GO is constructed. Based on this state transition chain, it is proven that any initial state of a search agent can reach the state corresponding to the GO with a positive probability through finite iterations, thereby ensuring that PSO escapes from LO with a probability of 1. In addition, it is deduced that increasing the distances between pbest and gbest, as well as the values of the learning factors $C_1, C_2$, enhances the speed at which PSO escapes from LO. Finally, numerical experiments are conducted to validate the theoretical findings.

The rest of this paper is organized as follows. Section 2 formulates the problem and establishes the theoretical framework. Section 3 presents the two core contributions of this work, along with their interpretations. Section 4 provides detailed proofs of the theoretical results. Section 5 validates the findings through numerical experiments. Finally, Section 6 concludes this work.

## 2. PROBLEM FORMULATION

This section formulates the problem for analyzing the ability of PSO to escape from LO. Section 2.1 introduces the solution update equation of PSO, Section 2.2 proposes the benchmark modelling, and Section 2.3 completes the problem formulation for analyzing PSO's ability to escape from LO.

### 2.1 Particle Swarm Optimization

PSO is designed to solve the following single-objective optimization problem:

$$\min_{\mathbf{lb} \leq \mathbf{x} \leq \mathbf{ub}} f(\mathbf{x}) \tag{2.1}$$

where **lb** and **ub** are the lower bound and upper bound vectors of decision variables, respectively, $f: \mathbb{R}^D \to \mathbb{R}$ denotes the objective function, and $D$ denotes the dimension of the optimization problem.

In PSO, the positions of search agents are updated as follows [6]:

$$v_{ij}(t+1) = \omega v_{ij}(t) + C_1 r_1 \left(pb_{ij}(t) - x_{ij}(t)\right) + C_2 r_2 \left(gb_j(t) - x_{ij}(t)\right) \tag{2.2}$$

$$x'_{ij}(t+1) = x_{ij}(t) + v_{ij}(t+1) \tag{2.3}$$

where $\mathbf{x}_i(t) = [x_{ij}(t)](i = 1, \ldots, N, j = 1, \ldots, D, t = 1, \ldots, T$, similarly throughout) and $\mathbf{v}_i(t) = [v_{ij}(t)]$ are the position and velocity vectors of search agent $i$ at iteration $t$, respectively. In this study, $N$ and $T$ represent the number of search agents and the maximum number of iterations, respectively. The vectors $\mathbf{pb}_i(t) = [pb_{ij}(t)]$ and $\mathbf{gb}(t) = [gb_j(t)]$ represent the historical best positions of agent $i$ and the entire population, respectively, and are referred to as *pbest* and *gbest*. $\omega$ denotes the *inertia weight*, and $C_1, C_2$ denote the *learning factors*. To date, scholars have proposed various mechanisms for selecting these three parameters [14-16,43-46]. In this study, the parameters $\omega, C_1$, and $C_2$ are considered constants that remain independent of the iterations, and $0 < \omega \leq 1, C_1, C_2 > 0$ [14,47]. The expression $\mathbf{x}'_i(t+1) = [x'_{ij}(t+1)]$ represents the updated position, which will be constrained into the box region to obtain the final position $\mathbf{x}_i(t+1) = [x_{ij}(t+1)]$ [48]:

$$x_{ij}(t+1) = B_{[lb_j, ub_j]_I}\left(x'_{ij}(t+1)\right) \tag{2.4}$$

where

$$B_{[p,q]_I}(x) = \begin{cases} p, x \leq p \\ x, p < x < q \\ q, x \geq q \end{cases} \tag{2.5}$$

It is important to note that in the existing literature, both vectors and closed intervals are represented using square brackets '[ ]'. To differentiate between the two in this paper, we use square brackets '[ ]' for vectors and square brackets with a subscript $I$, denoted as '[ ]$_I$', for closed intervals. For instance, $[lb_j, ub_j]_I$ in Eq. (2.4) and $[p,q]_I$ in Eq. (2.5) both represent closed intervals.

**2.2 Benchmark Modelling**

This section presents the benchmark modelling used to analyze the ability of PSO to escape from LO.

As previously mentioned, in PSO, each agent updates its position based on its pbest and the gbest. This results in thorough exploration of the spaces surrounding pbest and gbest, leading to the convergence of each agent's pbest towards a LO of the optimization problem (or, optimistically, the GO) within several iterations, while the gbest represents the best position among these pbests. Thereafter, although the positions of the agents continue to update, the positions of their pbests (and gbest) change little unless a better position is found, which corresponds to either the GO or a better LO. It is important to note that escaping from a LO does not necessarily mean that the algorithm has achieved the GO (rather, it may have converged to a better LO). Nevertheless, the ability to escape from LO is clearly a necessary condition for ensuring the global convergence ability of PSO.

As a result, this study analyzes the behavior of search agents after their pbests and the gbest have converged to specific local optima. We employ the well-known *stagnation assumption*: each agent's pbest and the population's gbest in PSO have stagnated at certain local optima [49-50]. Based on this assumption, if an agent discovers a better LO or GO during iterations, it can be asserted that the agent has successfully escaped from the LO. If any agent in the population escapes from the LO, it follows that the entire population has also escaped. This implies that, prior to escaping from the LO, the pbest and gbest can be considered constants that are independent of the iterations. Thus, in the subsequent analysis, $\mathbf{pb}_i(t)$ and $\mathbf{gb}(t)$ in Eq. (2.2) are reduced to $\mathbf{pb}_i$ and $\mathbf{gb}$, respectively. Furthermore, each $\mathbf{pb}_i$ or $\mathbf{gb}$ can be regarded as a specific LO.

*Case 1 (validity of the proposed benchmark model)*

This case provides a specific example to validate the effectiveness of the proposed benchmark model. In this instance, the PSO is utilized to optimize the well-known Rastrigin's Function, defined as $f(\mathbf{x}) = \sum_{i=1}^{D}(x_i^2 - 10\cos 2\pi x_i + 10)$ [51]. It is known that the Rastrigin's Function achieves its GO at $\mathbf{x} = \mathbf{0}$ and has multiple local optima, each LO approximately distributed around integer points [51]. For this analysis, we set $D = 2, N = 5, T = 200, \mathbf{lb} = [-5, -5]$, and $\mathbf{ub} = [5,5]$. Figures 1 (a) and (b) plot the first and second dimensions of each agent's pbest (i.e., $pb_{i1}(t)$ and $pb_{i2}(t)$), respectively. As shown in Fig. 1, after several iterations, each $pb_{ij}(t)$ approaches an integer, indicating that the pbest of each agent is located near a LO. Subsequently, each $pb_{ij}(t)$ either remains relatively unchanged over several iterations (indicating that it has become trapped in this LO) or shifts to another integer vicinity after several iterations (indicating that it has escaped from the LO and identified a better LO or GO). Finally, all pbests converge near the GO, demonstrating that the PSO successfully resolved the global optimization problem. Since the gbest is the best pbest, the conclusions drawn above also apply to the gbest. In summary, the effectiveness of the proposed benchmark model is validated.

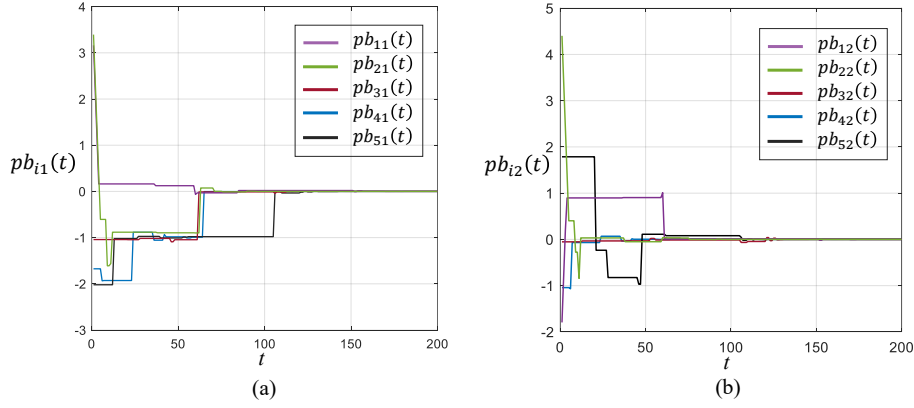

Fig. 1. The pbest of each agent. (a) $pb_{i1}(t)$, (b) $pb_{i2}(t)$.

The proposed benchmark model can be significantly simplified in two ways. First, the analysis of the population escaping from LO can be reduced to an examination of the behavior of a single agent. Thus, in the subsequent sections, we consider only one agent by setting $N = 1$. Second, since each dimension of the agent is updated independently, an agent escapes from the LO only when all of its dimensions have escaped. This indicates that the analysis of the agent escaping from LO can be further distilled to the analysis of the agent's behavior in a single dimension. Therefore, we consider only one dimension by setting $D = 1$. Consequently, the subscripts $i$ representing each agent and $j$ representing each dimension can be omitted from the variables that denote position and velocity in PSO. As a result, the previously defined optimization problem and the position update equation of PSO can be expressed in the following simplified form:

$$\min_{lb \leq x \leq ub} f(x) \tag{2.6}$$

$$v(t+1) = \omega v(t) + C_1 r_1 (pb - x(t)) + C_2 r_2 (gb - x(t)) \tag{2.7}$$

$$x'(t+1) = x(t) + v(t+1) \tag{2.8}$$

$$x(t+1) = B_{[lb,ub]_I}(x'(t+1)) \tag{2.9}$$

where $pb$ and $gb$ represent the stagnated positions of pbest and gbest, respectively, with each corresponding to a specific LO.

The subsequent proofs and analyses will be conducted on the basis of Eqs. (2.6)-(2.9). In addition, recognizing the symmetry between the terms $C_1 r_1 (pb - x(t))$ and $C_2 r_2 (gb - x(t))$, it is assumed that $pb \leq gb$ holds for the rest of this paper.

**2.3 Problem modelling for analyzing the ability of PSO to escape from LO**

Building on the previous benchmark modelling, this section establishes a model for analyzing the ability of PSO to escape from LO. This modeling process consists of three steps: First, we characterize the position and velocity of the agent as a Markov chain; next, we define the set of states corresponding to agents escaping from LO; finally, we provide a quantitative characterization of the agent's ability to escape from LO and outline the problems to be addressed in the subsequent sections.

Initially, the velocity and position of the agent are modelled using stochastic process theory. Although the ability to escape from LO pertains solely to the agent's position, the position of the agent at iteration $t + 1$ is influenced by both its position and velocity at iteration $t$ (as described in Eqs. (2.7)-(2.9)). Therefore, it is essential to jointly model the position and velocity of the agent. Clearly, during the iteration process, the position and velocity of the agent constitute a stochastic process. Within this framework, we define the state vector of the agent as follows:

$$\xi(t) = [x(t), v(t)] \in S_p \tag{2.10}$$

where $\xi(t)$ represents the state vector of the agent. $x(t)$ and $v(t)$ represent the position and velocity of the agent, respectively. $S_p$ denotes the search space of the agent, which can be expressed as follows:

$$S_p = [lb, ub]_I \times \mathbb{R} \tag{2.11}$$

where $'\times'$ denotes the Cartesian product of two sets: the first set $[lb, ub]_I$ represents the value range of $x(t)$, while the second set $\mathbb{R}$ represents the value range of $v(t)$.

According to Eqs. (2.7)-(2.9), the value of $\xi(t + 1)$ depends solely on $\xi(t)$. Thus, we have the following proposition:

*Proposition 1*

The state vector $\xi(t)$ constitutes a continuous state Markov Chain within the state space $S_p$.

Next, we present a rigorous formulation for the agent that successfully escapes from the LO. Let $S_g$ denote the set of states corresponding to the agent that has escaped from the LO:

$$S_g = R_g \times \mathbb{R} \subset S_p \tag{2.12}$$

Here, similar to the expression for $S_p$ in Eq. (2.11), $S_g$ is defined as the Cartesian product of two sets, where $R_g$ and $\mathbb{R}$ correspond to the positions and velocities of all agents that have escaped from the LO, respectively. Since escaping LO imposes a condition only on the agent's position, the right-hand side of the Cartesian product is $\mathbb{R}$. The following Case 2 provides an illustration for the meaning and construction method of $R_g$ in this work.

*Case 2 (Illustration of $R_g$)*

This case provides an example to illustrate the expression of $R_g$. We still consider the 1-D Rastrigin's Function $f(x) = x^2 - 10 \cos 2\pi x + 10$ [51] within the interval $[lb, ub]_I = [-2.5, 2.5]_I$. Obviously, $f(x)$ exhibits four local optima: $x_{l1} \approx -1.990, x_{l2} \approx 0.995, x_{l3} \approx 0.995, x_{l4} \approx 1.990$, along with one GO $x_g = 0$ within $[lb, ub]_I$. Supposing that the agent stagnates near $x_{l4}$, where $pb = x_{l4}$, the agent will not escape from the LO (i.e. $x_{l4}$) until a solution satisfying $f(x) \leq f(pb)$ is found. The inequality $f(x) \leq f(pb)$ defines the region $x \in R_2 \cup R_3 \cup R_5$, where $R_i$ represents the neighborhood of the LO $x_{li}(i = 1, \ldots, 4)$, and $R_5$ is the neighborhood of the GO. Consequently, in this case, $R_g = R_2 \cup R_3 \cup R_5$, which is the union of the neighborhoods of all better local optima and the GO. This process is illustrated in Fig. 2.

In the following discussion, to simplify the analysis, $R_g$ *is regarded as the neighborhood of the GO*, allowing $R_g$ to be represented as an interval, which can be expressed as follows:

$$R_g = [l_g, u_g]_I \tag{2.13}$$

On the one hand, this interval is clearly a subset of the solution set for $f(x) \leq f(pb)$, indicating that if the agent is able to locate a position within this interval, it surely bears the ability to escape from the LO. On the other hand, when the agent stagnates near the best LO, the solution set for $f(x) \leq f(pb)$ coincides with the neighborhood of the GO. In summary, this simplification does not compromise the validity of the proposed model.

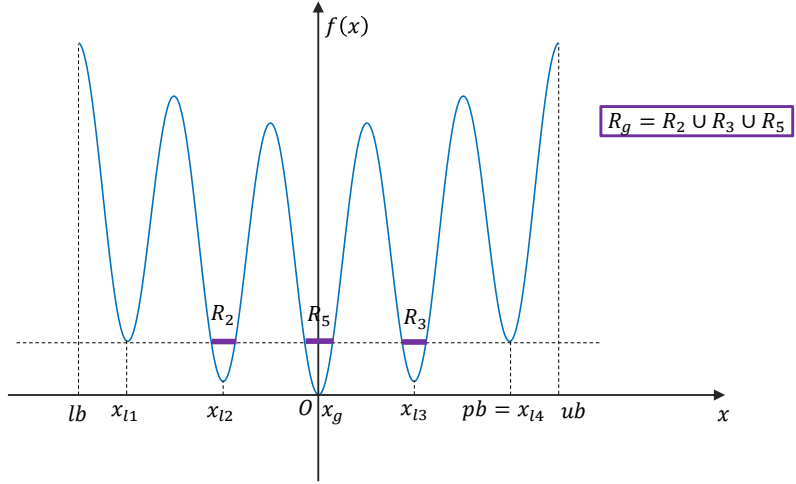

Fig. 2. $R_g$ in 1-D Rastrigin's Function

Finally, we present a quantitative characterization of the agent's ability to escape from LO.

To achieve this, we first introduce the concept of first hitting time. Let $S \subset S_p$ represent a set of states for the agent, and let $T(S|\xi(t_0))$ denote the first hitting time of the agent to the set $S$. It is evident that $T(S|\xi(t_0))$ satisfies the following formula:

$$\mathbb{P}\{T(S|\xi(t_0)) = t\} = \mathbb{P}\{\xi(t_0 + t) \in S, \xi(t_0 + \tau) \notin S, \forall \tau = 0, \dots, t-1|\xi(t_0)\} \quad (2.14)$$

and for fixed values of $\omega, C_1$, and $C_2$ that are independent of the iterations, the value of $T(S|\xi(t_0))$ is independent of $t_0$. Here, the blackboard bold symbol $\mathbb{P}$ in Eq. (2.14) represents probability, as used consistently throughout this paper.

In order to quantitatively analyze the agent's ability to escape from LO, we focus on the probability that the search agent successfully escapes from the LO within a finite number of iterations, denoted as $p_e$. Based on the aforementioned concept of first hitting time, $p_e$ can be expressed as follows, where $S_g$ is defined in Eqs. (2.12)-(2.13):

$$p_e = \lim_{t \to \infty} \mathbb{P}\{T(S_g|\xi(0)) \leq t\} \quad (2.15)$$

Clearly, the value of $p_e$ directly measures the ability of PSO to escape from the LO. If $p_e = 0$, it indicates that PSO has no ability to escape from the LO; if $0 < p_e < 1$, PSO possesses some ability to escape from the LO, but there remains a probability that it may never escape, regardless of the number of iterations. Conversely, if $p_e = 1$, PSO is guaranteed to escape from the LO within a finite number of iterations, underscoring its strong capability to escape from LO. In addition, calculating $p_e$ relies heavily on analyzing the behavior of the search agent. Thus, once the value of $p_e$ is determined, the analysis of the potential behaviors of the agent prior to escaping from the LO is inherently concluded.

In summary, the following two questions will be addressed in the subsequent sections:

(1) What is the value of $p_e$?

(2) If $p_e > 0$, what might be the behavior of the search agent before escaping from the LO?

## 3. MAIN THEORETICAL RESULTS

This section addresses solutions to the two questions posed at the end of Section 2. Specifically, the core theoretical results of this study are Proposition 2 and Proposition 3, with the latter being directly derivable from the former through straightforward probabilistic analysis. We present these two conclusions as follows (where $S_g$ is defined in Eqs. (2.12)-(2.13), and $T\left(S_g|\boldsymbol{\xi}(t_0)\right)$ represents the first hitting time from $\boldsymbol{\xi}(t_0)$ to $S_g$, as defined in Eq. (2.14)):

*Proposition 2*

When $\omega = 1$ and $pb \neq gb$, there exist $p_{e0} > 0$ and $t_{e0} \in \mathbb{N}^+$ such that for any given state vector $\boldsymbol{\xi}(t_0)$ at any initial iteration $t_0$, it holds that

$$\mathbb{P}\left\{T\left(S_g|\boldsymbol{\xi}(t_0)\right) \leq t_{e0}\right\} \geq p_{e0} \tag{3.1}$$

*Proposition 3*

When $\omega = 1$ and $pb \neq gb$, then for any initial state $\boldsymbol{\xi}(0)$, $p_e = 1$ holds.

Under the conditions that $\omega = 1$ and $pb \neq gb$, Proposition 2 implies that for any initial state $\boldsymbol{\xi}(t_0)$, the search agent can find a state within $S_g$ in at most $t_{e0}$ iterations, with a probability of at least $p_{e0}$. Building on Proposition 2, Proposition 3 states that the search agent can escape from LO in a finite number of iterations with a probability of 1, regardless of its initial state. The proof of Proposition 3 is provided in Appendix A.1 (based on Proposition 2), while the proof of Proposition 2 will be elaborated in Section 4.

We provide several clarifications regarding these two core theoretical results:

(a) A key condition for the validity of Proposition 2 and Proposition 3 is that the inertia weight $\omega = 1$ in PSO. In Section 4, we will elaborate on the role of $\omega = 1$ in the proof process of Proposition 2 and analyze PSO's ability to escape from LO ability when $\omega < 1$.

(b) Another assumption necessary for the validity of Proposition 2 and Proposition 3 is that $pb \neq gb$. It is important to note that under the stagnation assumption given in Section 2.2, $pb$ and $gb$ represent two stagnated local optima of the optimization problem. As a result, Proposition 3 indicates that, when *the agent's pbest and the population's gbest are stagnated around different local optima, the agent can escape from the LO in finite iterations with a probability of 1*; thus, its ability to escape from LO is well guaranteed. Section 4 will provide a comprehensive discussion on the significance of $pb \neq gb$ in the proof process of Proposition 2 and will also analyze PSO's ability to escape from LO when $pb = gb$.

(c) As discussed in Section 1, although for fixed agent's velocity and position at iteration *t*, the agent position's distribution may not cover the entire feasible region at iteration $t + 1$, the distribution at iterations $t + 2, t + 3$, and beyond remains uncertain, indicating that the global searching capability of the agent in subsequent iterations is also uncertain. Interpreting from the *perspective of probability distribution*, Proposition 3 suggests that regardless of the initial state of the agent, its search range can cover the entire feasible region after several iterations, thereby ensuring that the probability of escaping from the LO in each iteration is greater than zero. Consequently, PSO can escape from LO in a finite number of iterations with a probability of 1. This assertion will be validated through numerical simulations in Section 5.

(d) Interpreting from the *perspective of Markov chain*, Proposition 2 asserts that there exists a

positive integer $t_{e0}$ such that for any initial state $\xi(t_0)$, the $t_{e0}$-step transition probability from $\xi(t_0)$ to $S_g$ is positive. In other words, $S_g$ is accessible from any possible initial state $\xi(t_0)$. To prove Proposition 2, it is sufficient to establish a state transition chain from any potential initial state set $S_0 = \{\xi(t_0)\}$ to $S_g$. Specifically, the goal is to find a series of state sets $S_t \subset S_p$ such that the state transition chain of $\xi(t): S_0 \to S_1 \to \cdots \to S_{t_{e0}-1} \to S_g$ exists. On this basis, by setting:

$$p_{e0} = \mathbb{P}\{\xi(t_0 + t_{e0}) \in S_g, \xi(t_0 + t) \in S_t, t = 1, \ldots, t_{e0} - 1 | \xi(t_0)\} > 0 \quad (3.2)$$

it holds that $\mathbb{P}\{T(S_g|\xi(t_0)) \leq t_{e0}\} \geq p_{e0} > 0$. In summary, $t_{e0}$ represents an upper bound on the first hitting time from $S_0$ to $S_g$, while $p_{e0}$ serves as a lower bound for the transition probability of the corresponding transition chain.

From this perspective, Section 4 will prove Proposition 2 by constructing a state transition chain that satisfies Eq. (3.2) and will analyze the potential behaviors of the agent prior to escaping from the LO during the proof process.

## 4. PROOF AND DISCUSSIONS OF PROPOSITION 2

This section will present a rigorous proof of Proposition 2 and analyze the potential behavior of the search agent prior to escaping from LO. It is important to note that throughout this section, unless otherwise stated, it is assumed that $\omega = 1$ and $pb \neq gb$ (i.e. $pb < gb$, as postulated in Section 2.2).

**4.1 Preliminaries: Transition Kernel Density**

To calculate the values of $p_{e0}$ and $t_{e0}$ in Eq. (3.1), it is imperative to examine the properties of the Markov chain formed by the agent's state vector $\xi(t)$. For the continuous state Markov chain $\xi(t)$, which lacks a state transition matrix like that of a discrete state Markov Chain, we cannot apply matrix analysis techniques in the same manner as with discrete state Markov Chains. Instead, we will utilize the transition kernel and transition kernel density to characterize the properties of $\xi(t)$. These properties will enhance the understanding of the solution update operator in PSO and will play a significant role in the calculations and analyses presented in subsequent sections.

We first provide the definitions of the transition kernel and transition kernel density for $\xi(t)$:

*Definition 1 (transition kernel and transition kernel density)*

For a specific state $\zeta \in S_p$ and a state set $S \subset S_p$, the transition kernel of $\xi(t)$ is defined as:
$$F(\zeta, S) = \mathbb{P}\{\xi(t+1) \in S | \xi(t) = \zeta\} \quad (4.1)$$
In addition, the transition kernel density $f(\zeta, \eta)$ is defined as a non-negative function such that:

$$F(\zeta, S) = \int_{\eta \in S} f(\zeta, \eta) \mathrm{d}\eta \quad (4.2)$$

It is evident that the state transition probability of a discrete-state Markov chain is the discretized version of the transition kernel from a continuous-state Markov chain. As the derivative of the transition kernel, the transition kernel density possesses a distinct probabilistic interpretation: for a fixed $\xi(t) = \zeta$, $f(\zeta, \eta)$ represents the joint Probability Density Function (joint PDF) of the random variables $x(t+1)$ and $v(t+1)$ (noting that $\eta$ is a two-dimension vector, with the first dimension corresponding to the agent's position and the second corresponding to the agent's velocity). In addition, since the value of $x(t + 1)$ is determined solely by $v(t + 1)$ and $x(t)$ (as referenced in Eqs. (2.8)-(2.9)), we can simplify the problem of determining the transition kernel density of $\xi(t)$ to finding the PDF of the random variable $v(t + 1)$ under the condition that $\xi(t) = \zeta$.

According to Eq. (2.7), for fixed $\xi(t)$, $v(t + 1)$ can be expressed as the sum of two random variables: $v(t) + C_1 r_1(pb - x(t))$ and $C_2 r_2(gb - x(t))$, both of which follow uniform distributions. Consequently, the PDF of $v(t + 1)$ can be derived as follows (see Appendix A.3 for the proof of Proposition 4):

*Proposition 4*

Given the state vector $\xi(t)$, the PDF of $v(t + 1)$ (denoted as $f_{v(t+1)}(v|\xi(t))$) can be expressed as follows:

$$f_{v(t+1)}(v|\xi(t)) = \begin{cases} \frac{h_f(v - v_{f_1})}{v_{f_2} - v_{f_1}}, & v_{f_1} \leq v < v_{f_2} \\ h_f, & v_{f_2} \leq v < v_{f_3} \\ -\frac{h_f(v - v_{f_4})}{v_{f_4} - v_{f_3}}, & v_{f_3} \leq v < v_{f_4} \\ 0, & \text{otherwise} \end{cases} \quad (4.3)$$

where $v_{f_i}(i = 1, \ldots, 4)$ and $h_f$ are parameters expressed by $pb, gb, v(t)$, and $x(t)$, as described in Appendix A.2.

According to Proposition 4, the curve of $f_{v(t+1)}(v|\xi(t))$ along with the $v$ axis form an isosceles trapezoid. The curve of $f_{v(t+1)}(v|\xi(t))$ is shown in Fig. 3.

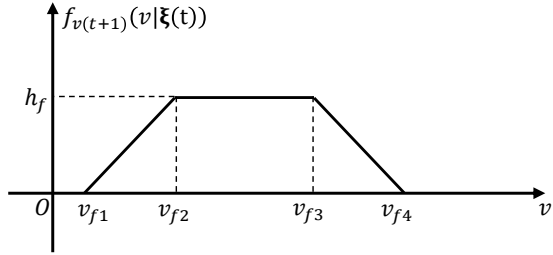

Fig. 3. Curve of $f_{v(t+1)}(v|\xi(t))$

Based on Proposition 4, several properties of the distribution of $v(t + 1)$ are given as follows (see Appendix A.4 for the proof of Corollary 1).

*Corollary 1*

Given $\xi(t)$, by defining the range of $v(t + 1)$ as $D_{v(t+1)}$, the following properties hold:

(a) $D_{v(t+1)} = [v_{f_1}, v_{f_4}]_I$, and $v(t) \in D_{v(t+1)}$.

(b) Set
$$d_0 = \min\{C_1, C_2, 1\}(gb - pb) \quad (4.4)$$
Then it follows that
$$|D_{v(t+1)}| \geq d_0 \quad (4.5)$$
(c) For any $\xi(t)$ and interval $D_1 \subset D_{v(t+1)}$, it holds that
$$\mathbb{P}\{v(t + 1) \in D_1\} \geq h(|D_1|) > 0 \quad (4.6)$$
where

$$h(x) = \min \left\{ \frac{1}{2\max\{C_1,C_2\}(ub-lb)} x, \frac{1}{2C_1 C_2 (ub-lb)^2} x^2 \right\} \tag{4.7}$$

Several remarks regarding Corollary 1 are as follows:

Corollary 1 (a) points out that the distribution range of $v(t+1)$ covers $v(t)$. Specifically, setting $r_1 = r_2 = 0$ in Eq. (2.7) directly yields $v(t+1) = v(t)$ under the condition that $\omega = 1$.

Corollary 1 (b) points out that the interval length of the distribution range of $v(t+1)$ is always no less than $d_0 = \min\{C_1, C_2, 1\}(gb - pb)$.

Corollary 1 (c) points out that, for any set $D_1$ with non-zero measure within the distribution range of $v(t+1)$, the probability that $v(t+1)$ lies within $D_1$ is positive. According to Eqs. (2.8)-(2.9), this corollary also applies to $x(t+1)$. This corollary will play a central role in subsequent probability analysis.

**4.2 Analysis of possible behavior of the search agent**

This section introduces two typical behaviors of the search agent, serving as the basis for construction of the state transition chain for $\xi(t)$ that satisfies Eq. (3.2) in Sections 4.3-4.5. On these bases, the explicit expressions for $t_{e0}$ and $p_{e0}$ will be provided in Section 4.6.

First, consider the following two simple cases:

*Case 3*

Set $\xi(0) = [1,1], R_g = [5.8, 6.2]_I, [lb, ub]_I = [0, 10]_I$. As previously noted, setting $r_1 = r_2 = 0$ in Eq. (2.7) yields $v(t+1) = v(t)$, implying that the agent maintains a constant rightward velocity $v = v(0) = 1$. Consequently, by continuously setting $r_1 = r_2 = 0$, $x(t)$ forms an arithmetic sequence with the first term $x(0) = 1$ and a common difference $v(0) = 1$, leading to $x(6) = 6 \in R_g$. Therefore, it is reasonable to assert that the probability of the agent escaping from LO after the 6$^{th}$ iteration is positive.

*Case 4*

Set $\xi(0) = [8,1], R_g = [5.8, 6.2]_I, [lb, ub]_I = [0, 10]_I$. In this case, continuously setting $r_1 = r_2 = 0$ yields $x(1) = x(0) + v(1) = x(0) + v(0) = 9$, and $x(2) = x(1) + v(2) = x(1) + v(0) = 10 = ub$. Subsequently, setting $r_1 = r_2 = 1$ in Eq. (2.7) gives us $v(t+1) = v(t) - C_1(ub - pb) - C_2(ub - gb)$, implying that the value of $v(t)$ will continuously decrease until it becomes negative. Afterwards, by setting $r_1 = r_2 = 0$ in Eq. (2.7), the agent consistently maintains a leftward velocity, and the search agent will continuously move towards left until it finds a position within $R_g$.

It should be noted that the probability that $r_1 = r_2 = 1$ or $r_1 = r_2 = 0$ is zero. Therefore, the probabilities of the scenarios described in above cases are also zero. However, the movement trends of the search agents in these two cases are entirely plausible.

From these two cases, it can be concluded that the behavior of the search agent prior to escaping from LO can be categorized into the following two typical types.

(1) *Inertial motion*: When $r_1$ and $r_2$ are close to 0, the sign of the agent's velocity remains unchanged, allowing it to continue moving in a specific direction until it reaches the upper or lower bound of the feasible region. It's important to note that this property holds when $\omega = 1$; otherwise, the agent's velocity in that direction will gradually decay, potentially causing it to lose the ability to continuously move towards the upper (or lower) boundary of the feasible region.

(2) *Oscillation*: When $r_1$ and $r_2$ are close to 1, if the agent is moving to the right and $x(t) >$

max{$pb, gb$}, the negative terms $C_1r_1(pb - x(t)) + C_2r_2(gb - x(t))$ will cause the agent's velocity to decrease until it becomes negative, prompting the agent to change direction and move left, thus completing a leftward oscillation. Conversely, if the agent is moving to the left and $x(t) <$ max{$pb, gb$}, symmetry suggests a probability that the agent will complete a rightward oscillation.

*Case 5*

This case illustrates a scenario where $\omega < 1$ and the agent fails to reach the upper bound $ub$ within a single rightward inertial motion. Suppose that $\boldsymbol{\xi}(0) = [0,1], [lb, ub]_I = [-2,5]_I, pb = -1, gb = -2$, and $\omega = 0.8$. In this situation, the agent is performing a single rightward inertial motion. Under these assumptions, it holds that $\forall t \in \mathbb{N}^+, v(t) > 0$, and $x(t) > x(t-1)$. Consequently, we can conclude that $\forall t \in \mathbb{N}^+, x(t) > x(0) > \max\{pb, gb\}$, indicating that $\forall t \in \mathbb{N}^+, v(t) = 0.8v(t-1) + C_1r_1(pb - x(t-1)) + C_2r_2(gb - x(t-1)) \leq 0.8v(t-1)$, thus it holds that $\forall t \in \mathbb{N}^+, x(t) = x(0) + \sum_{\tau=1}^{t} v(\tau) \leq \sum_{\tau=1}^{t} 0.8^\tau < 4 < ub$. In summary, when $\omega < 1$, the agent may fail to reach the upper bound $ub$ within a single rightward inertial motion.

In Case 3, the agent escapes from the LO after a single (rightward) inertial motion. In Case 4, the agent escapes from the LO after a sequence of one (rightward) inertial motion, one (leftward) oscillation, and one (leftward) inertial motion. It is important to note that prior to escaping from LO, the agent is likely to undergo multiple inertial motions and oscillations. However, for the sake of simplifying the analysis, this section will only consider the agent's behavior as described in Case 4. This means that in subsequent sections, the state transition chain from $S_0 = \{\boldsymbol{\xi}(t_0)\}$ to $S_g$ will be constructed through the following three steps:

Step 1: The agent completes a *rightward inertia motion* from the state set $S_0$ to $S_a$, where
$$S_a = \{\boldsymbol{\xi}(t) | x(t) = ub, v(t) > 0\} \tag{4.8}$$
Obviously, upon reaching $S_a$, the position of the agent attains the upper bound $ub$, thereby completing the rightward inertia motion.

Step 2: The agent undergoes a *leftward oscillation motion* from the state set $S_a$ to $S_b$, where
$$S_b = \{\boldsymbol{\xi}(t) | x(t-1) = ub, v(t) < 0\} \tag{4.9}$$
During this oscillation motion, the position of the agent remains unchanged (equal to $ub$), while its velocity continuously decreases from a positive value to a negative value. After reaching $S_b$, the position of the agent changes ($x(t) = x(t-1) + v(t) < ub$), thus completing the leftward oscillation motion.

Step 3: The agent completes a *leftward inertia motion* from the state set $S_b$ to $S_g$, where $S_g$ is defined in Eqs. (2.12)-(2.13), representing the set of states corresponding to the agent that has escaped from the LO.

The construction methods for the state transition chains from Step 1 to Step 3 will be presented in Sections 4.3 to 4.5, respectively. Building on these state transition chain results, we will derive the expressions for the upper bound of the first hitting time and the lower bound of corresponding state transition probability. By combining the construction results of the state transition chains from Step 1 to Step 3, we will complete the final proof of Proposition 2, as outlined in Section 4.6.

### 4.3 Step 1: rightward inertial motion

This section aims to construct a state transition chain from $S_0$ to $S_a$, corresponding to a single rightward inertia motion. On this basis, an upper bound $t_{0a}$ for the first hitting time from $S_0$ to $S_a$ and a lower bound $p_{0a}$ for the corresponding transition probability can be determined.

Initially, in order to maintain the rightward inertia motion, the agent should possess a certain rightward initial velocity. However, this may not always hold due to the arbitrariness of the agent's initial state. Nevertheless, it can be proved that (see Appendix A.5 for the proof of Lemma 1):

*Lemma 1*

For any state $\xi(t)(t \geq 0)$, the following inequality holds:

$$\mathbb{P}\left\{v(t+1) \geq \tfrac{1}{4}d_0 \text{ or } v(t+1) \leq -\tfrac{1}{4}d_0\right\} \geq h\left(\tfrac{1}{4}d_0\right) \tag{4.10}$$

where $d_0$ is defined in Eq. (4.4), and $h(\cdot)$ is defined in Eq. (4.7).

According to Lemma 1, regardless of the value of the agent's state $\xi(t)$ at iteration $t$, there is a certain probability (no less than $h\left(\tfrac{1}{4}d_0\right)$) at iteration $t+1$ that the agent will exhibit either a rightward velocity (not less than $\tfrac{1}{4}d_0$) or a leftward velocity (not less than $\tfrac{1}{4}d_0$). Therefore, by symmetry, for the initial state $\xi(t_0)(t_0 \geq 1)$, we can assume that $\mathbb{P}\left\{v(t_0) \geq \tfrac{1}{4}d_0\right\} \geq h\left(\tfrac{1}{4}d_0\right)$. This implies that, starting from iteration $t_0$, the agent has a certain probability of achieving a rightward velocity, thereby facilitating a rightward inertial motion. In summary, in subsequent discussions, we can define the initial state set as:

$$S_0 = \left\{\xi(t) | v(t) \geq \tfrac{1}{4}d_0\right\} \tag{4.11}$$

Next, for a given initial state $\xi(t_0) = [x(t_0), v(t_0)] \in S_0$, we present a specific method for constructing the state transition chain. We denote the state transition chain as $\{\xi(t_0)\} \to S_1^1 \to S_2^1 \to \cdots S_{t_{0a}-1}^1 \to S_a$, where

$$S_t^1 = R_t^1 \times \mathbb{R}, t = 1, \ldots, t_{0a} - 1 \tag{4.12}$$

Here, $R_t^1$ represent a series of increasing intervals corresponding to the range of values for the agent's position within $[lb, ub]_I$ (indicating that $\forall t = 1, \ldots, t_{0a} - 2, \inf R_{t+1}^1 > \sup R_t^1$ and $\sup R_{t_{0a}-1}^1 < ub$). The superscripts $'1'$ on $S_t^1$ and $R_t^1$ indicate their association with Step 1. This construction implies that we *will not analyze the agent's velocity but will instead focus solely on the range of values for the agent's position*. We also define $R_0^1 = \{x(t_0)\}, R_{-1}^1 = \{x(t_0) - v(t_0)\}$, and an interval $R_{t_{0a}}^1$ satisfying $\inf R_{t_{0a}}^1 \geq ub$. Fig. 4 provides a schematic diagram of this series of intervals $R_t^1(t = -1, \ldots, t_{0a})$.

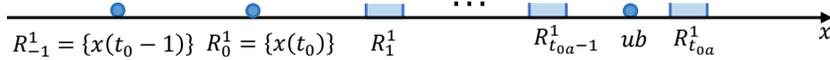

$R_{-1}^1 = \{x(t_0 - 1)\}$  $R_0^1 = \{x(t_0)\}$  $R_1^1$  $R_{t_{0a}-1}^1$  $ub$  $R_{t_{0a}}^1$

Fig. 4 Schematic diagram of $R_t^1(t = 1, \ldots, t_{0a})$.

Furthermore, we assume $x(t + t_0) \in R_t^1, t = -1, \ldots, t_{0a}$. Obviously, for the given $\xi(t_0) \in S_0$, when neglecting boundary constraints, both $x(t_0)$ and $x(t_0 - 1) = x(t_0) - v(t_0)$ are fixed, corresponding to the sets $R_0^1$ and $R_{-1}^1$. Subsequently, without considering boundary constraints, since $v(t + t_0) = x(t + t_0) - x(t + t_0 - 1)$, and given that $\xi(t)$ forms a Markov chain, the distribution range of $x(t + t_0 + 1)$ can be uniquely determined by the values of $x(t + t_0 - 1)$ and $x(t + t_0)$. If it can be proven that $\forall t \geq 0, \forall x(t + t_0 - 1) \in R_{t-1}^1, \forall x(t + t_0) \in R_t^1$, then the distribution range of $x(t + t_0 + 1)$ will always cover $R_{t+1}^1$. According to Corollary 1(c), it holds that $\forall t \geq 0$,

$$\inf_{x(t+t_0-1)\in R^1_{t-1}, x(t+t_0)\in R^1_t} \mathbb{P}\{x(t+t_0+1) \in R^1_{t+1} | x(t+t_0-1), x(t+t_0)\} \geq h(|R^1_{t+1}|) > 0$$

(4.13)

where $h(\cdot)$ is defined in Eq. (4.7). Finally, once the agent's position reaches $R^1_{t_{0a}}$, it is constrained to $ub$ and successfully attains the state set $S_a$. This directly implies that the state transition chain $\{\xi(t_0)\} \to S^1_1 \to S^1_2 \to \cdots S^1_{t_{0a}-1} \to S_a$ exists. Therefore, the first hitting time from $\{\xi(t_0)\}$ to $S_a$ has an upper bound $t_{0a}$, and the corresponding transition probability has a lower bound:

$$Pr \geq \prod_{t=0}^{t_{0a}-1} \inf_{x(t+t_0-1)\in R^1_{t-1}, x(t+t_0)\in R^1_t} \mathbb{P}\{x(t+t_0+1) \in R^1_{t+1} | x(t+t_0-1), x(t+t_0)\}$$

$$= \prod_{t=1}^{t_{0a}} h(|R^1_t|)$$

(4.14)

Based on the analysis above, for any given initial state $\xi(t_0) \in S_0$, we present the explicit expression for the state transition $\{\xi(t_0)\} \to \cdots \to S_a$ along with the corresponding values of $t_{0a}$ and $p_{0a}$ in the following Lemma 2, which is rigorously proven in Appendix A.6-Appendix A.8.

It is important to note that in the following Lemma 2, $t_1, t_2, a_t, \delta_a, b_t, \delta_b$ are temporary coefficients introduced to facilitate the description of the state transition chain. The definitions of these coefficients are applicable only to Lemma 2 and its proof. If these coefficients appear in subsequent lemmas or propositions, their meanings may not necessarily correspond to the definitions in Lemma 2. The reason for consistently using these coefficient names is that distinguishing the names of numerous minor coefficients in each lemma would result in overly complex notation, potentially hindering readability. Nevertheless, the naming of the state transition chains in Steps 1 to 3 will be distinct (denoted as $S^1_t, S^2_t, S^3_t$, respectively, as described in the explanations following Eq. (4.12)) to enhance clarity.

*Lemma 2*

When $\omega = 1$ and $pb \neq gb$, the following propositions hold:

(a) For a given $\lambda \in (0,1)$, if an initial state vector $\xi(t_0) \in S_0$ satisfies $x(t_0) \geq \lambda pb + (1-\lambda)gb$, then there exists a positive integer $t_1 \leq \left\lceil \frac{8(ub-lb)}{d_0} \right\rceil - 1$, such that the state transition chain $\{\xi(t_0)\} \to S^1_1 \to S^1_2 \to \cdots S^1_{t_1-1} \to S_a$ is feasible from $\{\xi(t_0)\}$ to $S_a$, where:

$$\begin{cases} S^1_t = R^1_t \times \mathbb{R}, t = 1, \dots, t_1 - 1 \\ R^1_t = [a_t - \delta_a, a_t]_I \\ a_t = x(t_0) + tv(t_0) - t(t-1)\delta_a \\ \delta_a = \min\left\{\frac{1-\lambda}{4}d_0, \frac{d_0^2}{128(ub-lb)}\right\} \end{cases}$$

(4.15)

(b) If an initial state vector $\xi(t_0) \in S_0$ satisfies $x(t_0) < \frac{1}{4}pb + \frac{3}{4}gb$, then there exists a positive integer $t_2 \leq \left\lceil \frac{4(ub-lb)}{d_0} \right\rceil$, such that the state transition chain $\{\xi(t_0)\} \to U^1_1 \to U^1_2 \to \cdots U^1_{t_2-1} \to U_m$ is feasible from $\{\xi(t_0)\}$ to $U_m$, where $U_m = \{\xi(t)|x(t) \geq \frac{1}{2}pb + \frac{1}{2}gb\}$, and:

$$\begin{cases} U_t^1 = W_t^1 \times \mathbb{R}, t = 1, \ldots, t_2 - 1 \\ W_t^1 = [b_t, b_t + \delta_b]_I \\ b_t = x(t_0) + tv(t_0) + t(t-1)\delta_b \\ \delta_b = \frac{1}{16} d_0 \end{cases} \quad (4.16)$$

(c) For any given state vector $\xi(t_0) \in S_0$ at any initial iteration $t_0$, it holds that:
$$\mathbb{P}\{T(S_a|\xi(t_0)) \leq t_{0a}\} \geq p_{0a} \quad (4.17)$$
where $S_a$ is defined in Eq. (4.8), and
$$\begin{cases} t_{0a} = 13 \left\lceil \frac{ub-lb}{d_0} \right\rceil \\ p_{0a} = h^{t_{0a}}\left(\frac{d_0^2}{128(ub-lb)}\right) \end{cases} \quad (4.18)$$
where $h(\cdot)$ is defined in Eq. (4.7).

We provide an explanation for the implications of Lemma 2. For any given $\lambda \in (0,1)$ and initial state vector $\xi(t_0) \in S_0$ satisfying $x(t_0) \geq \lambda pb + (1-\lambda)gb$, Lemma 2(a) outlines a method to construct a state transition chain from $\{\xi(t_0)\}$ to $S_a$. If the initial state $\xi(t_0)$ satisfies $x(t_0) \geq \frac{1}{4}pb + \frac{3}{4}gb$, a feasible state transition chain is provided by Lemma 2(a) (with $\lambda = \frac{1}{4}$).

Conversely, for the case where $x(t_0) < \frac{1}{4}pb + \frac{3}{4}gb$, the construction of the state transition chain occurs in the following two steps:

(1) Lemma 2(b) outlines the method for constructing a state transition chain from $\{\xi(t_0)\}$ to $U_m$, where $U_m = \{\xi(t)|x(t) \geq \frac{1}{2}pb + \frac{1}{2}gb\}$ represents the set of agents positioned at or above $\frac{1}{2}pb + \frac{1}{2}gb$.

(2) Considering the agent $\xi(t)$ within the state transition chain $\{\xi(t_0)\} \to U_1^1 \to U_2^1 \to \cdots U_{t_2-1}^1 \to U_m$ (as described in Lemma 2(b)), since the agent reaches $U_m$ at iteration $t_2 + t_0$, by the definition of $U_m$, it holds that $x(t_2 + t_0) > \frac{1}{2}(pb + gb)$. If we can further demonstrate that $\inf U_m - \sup U_{t_2-1}^1 \geq \frac{1}{4} d_0$, then it holds that $v(t_2 + t_0) = x(t_2 + t_0) - x(t_2 + t_0 - 1) \geq \frac{1}{4} d_0$, indicating that $\xi(t_2 + t_0) \in S_0$. Consequently, according to Lemma 2(a), we can construct a transition chain from $\{\xi(t_2 + t_0)\}$ (where $\xi(t_2 + t_0) \in U_m$) to $S_a$ commencing from iteration $t_2 + t_0$ (with $\lambda = \frac{1}{2}$ in Lemma 2(a)).

In conclusion, by combining the above two cases, we can derive an upper bound $t_{0a}$ for the first hitting time and a lower bound $p_{0a}$ for the corresponding first transition probability, as stated in Lemma 2(c). The rationale for treating the two cases $x(t_0) \geq \frac{1}{4}pb + \frac{3}{4}gb$ and $x(t_0) < \frac{1}{4}pb + \frac{3}{4}gb$ separately lines in the differing construction methods for the state transition chain required in each scenario. Detailed proofs of Lemma 2(a) to Lemma 2(c) can be found in Appendix A.6-Appendix A.8.

## 4.4 Step 2: Leftward Oscillation

This section aims to construct a state transition chain from $S_a$ to $S_b$, corresponding to a single leftward oscillation. On this basis, an upper bound $t_{ab}$ for the first hitting time from $S_a$ to $S_b$ and a lower bound $p_{ab}$ for the corresponding transition probability can be obtained.

Given the initial state vector $\xi(t_0) \in S_a$, as previously discussed, during the agent's leftward oscillation, its position remains constant at $ub$, while its velocity gradually decreases until it becomes negative. This implies that, for the agent to complete a leftward oscillation, the initial rightward velocity $v(t_0)$ must not be excessively high.

First, we provide an estimate for the upper bound of the agent's initial velocity $v(t_0)$, which is rigorously proven in Appendix A.9:

*Lemma 3*

$\forall t \geq 0$, the agent's velocity $v(t)$ has an upper bound $v_u$, where
$$v_u = (C_1 + C_2 + 1)(ub - lb) \tag{4.19}$$

Next, we present the specific method for constructing the state transition chain. Given the initial state vector $\xi(t_0) \in S_a$, since the position of the agent remains constant throughout the leftward oscillation (is equal to $ub$), we can denote the state transition chain as $\{\xi(t_0)\} \to S_1^2 \to S_2^2 \to \cdots S_{t_{ab}-1}^2 \to S_b$, where
$$S_t^2 = \{ub\} \times V_t^2, t = 1, \ldots, t_{ab} - 1 \tag{4.20}$$

Here, $V_t^2$ represents a series of decreasing intervals corresponding to the range of values for the agent's velocity (indicating that $\forall t = 1, \ldots, t_{ab} - 2, \sup V_{t+1}^1 < \inf V_t^1$ and $\sup V_1^1 < v(t_0)$). The superscripts '2' on $S_t^2$ and $V_t^2$ indicate their association with Step 2. The set $\{ub\}$ corresponds to the value of the agent's position. We also define $V_0^2 = \{v(t_0)\}$ and an interval $V_{t_{ab}}^2$ such that $V_t^2(t = 1, \ldots, t_{ab})$ forms a series of decreasing intervals, where $\sup V_{t_{ab}}^2 \leq 0$. Fig. 5 provides a schematic diagram of this series of intervals $V_t^2(t = 0, \ldots, t_{ab})$.

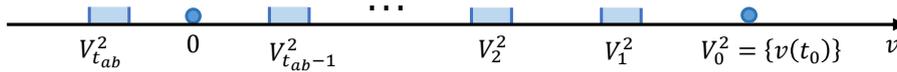

Fig. 5 Schematic diagram of $V_t^2(t = 0, \ldots, t_{ab})$.

Furthermore, we assume $v(t + t_0) \in V_t^2, t = 0, \ldots, t_{ab}$. It is evident that $\forall t = 0, \ldots, t_{ab} - 1, x(t + t_0) = ub$ holds; thus, the distribution range of $v(t + t_0 + 1)$ can be uniquely determined by the value of $v(t + t_0)$. Consequently, if we can demonstrate that $\forall t = 0, \ldots, t_{ab} - 1, \forall v(t + t_0) \in V_t^2$, the distribution range of $v(t + t_0 + 1)$ always covers $V_{t+1}^2$, then according to Corollary 1(c), it holds that $\forall t \geq 0$,
$$\inf_{v(t+t_0) \in V_t^2} \mathbb{P}\{v(t + t_0 + 1) \in V_{t+1}^2 | v(t + t_0)\} \geq h(|V_{t+1}^2|) > 0 \tag{4.21}$$

where $h(\cdot)$ is defined in Eq. (4.7). Finally, when the agent's velocity reaches $V_{t_{ab}}^2$, the agent gains a leftward velocity and successfully reaches the state set $S_b$. This indicates that the state transition chain $\{\xi(t_0)\} \to S_1^2 \to S_2^2 \to \cdots S_{t_{ab}-1}^2 \to S_b$ exists. Therefore, the first hitting time from $\{\xi(t_0)\}$ to $S_b$ is bounded above by $t_{ab}$, and the corresponding transition probability has a lower bound:
$$Pr \geq \prod_{t=0}^{t_{oa}-1} \inf_{v(t+t_0) \in V_t^2} \mathbb{P}\{v(t + t_0 + 1) \in V_{t+1}^2 | v(t + t_0)\} = \prod_{t=1}^{t_{oa}} h(|V_t^2|) \tag{4.22}$$

Based on the analysis above, for a given initial state vector $\xi(t_0) \in S_a$, we present the explicit expression for the state transition $\{\xi(t_0)\} \to \cdots \to S_b$ and the corresponding values of $t_{ab}$ and $p_{ab}$ in the following Lemma 4, which is rigorously proven in Appendix A.10:

*Lemma 4*

For two given real parameters $\lambda, \mu$ satisfying $0 < \lambda < \mu \leq \frac{1}{2}$, define a subset of $S_b$ as:

$$S_b(\lambda, \mu) = \{\xi(t) | x(t-1) = ub, -\mu d_0 < v(t) < -\lambda d_0\} \tag{4.23}$$

When $\omega = 1$ and $pb \neq gb$, given an initial state vector $\xi(t_0) \in S_a$, the following propositions hold:

(a) The transition chain $\{\xi(t_0)\} \to S_1^2 \to S_2^2 \to \cdots S_{\left\lceil\frac{3v(t_0)}{2d_0}\right\rceil}^2 \to S_b(\lambda, \mu)$ is feasible from $\{\xi(t_0)\}$ to $S_b(\lambda, \mu)$, where:

$$\begin{cases} S_t^2 = \{ub\} \times V_t^2, t = 1, \ldots, \left\lceil\frac{3v(t_0)}{2d_0}\right\rceil \\ V_t^2 = \left[v(t_0) - t\frac{v(t_0)}{\left\lceil\frac{3v(t_0)}{2d_0}\right\rceil}, v(t_0) - \left(t - \frac{1}{2}\right)\frac{v(t_0)}{\left\lceil\frac{3v(t_0)}{2d_0}\right\rceil}\right]_I \end{cases} \tag{4.24}$$

(b) It holds that

$$\mathbb{P}\{T(S_b(\lambda, \mu) | \xi(t_0)) \leq t_{ab}\} \geq p_{ab} \tag{4.25}$$

where

$$\begin{cases} t_{ab} = 2\left\lceil\frac{v_u}{d_0}\right\rceil \\ p_{ab} = h((\lambda - \mu)d_0)h^{t_{ab}}\left(\frac{1}{4}d_0\right) \end{cases} \tag{4.26}$$

We clarify the implications of Lemma 4. For $\xi(t_0) \in S_a$, rather than directly presenting a state transition chain from $\{\xi(t_0)\}$ to $S_b$, Lemma 4 presents a state transition chain from $\{\xi(t_0)\}$ to $S_b(\lambda, \mu)$, where $S_b(\lambda, \mu)$ is a subset of $S_b$, with additional constraints imposed on the agent's velocity.

The reason for defining the terminal state set of the left oscillation as $S_b(\lambda, \mu)$ rather than simply $S_b$ is that, in Section 4.5 (i.e. Step 3), when constructing the state transition chain from $S_b$ to $S_g$, the agent in $S_b$ must have an appropriate initial leftward velocity. Specifically, if this leftward velocity is too small, the agent may lose the capacity for leftward inertial motion; conversely, if the leftward velocity is too large, the agent may overshoot, complicating the process of precisely finding a position within the neighborhood of the GO (i.e. $R_g$, as discussed in Case 2). From a technical perspective (see Section 4.5), in Step 3, the different numerical relationships among $pb, gb$ and $R_g$ impose varying requirements on the initial leftward velocity needed to construct a feasible state transition chain. Therefore, we define $S_b(\lambda, \mu)$, a subset of $S_b$ with two undetermined coefficients $\lambda$ and $\mu$, as the terminal point of the state transition chain. This allows us flexibility in selecting different coefficients $\lambda$ and $\mu$ for $S_b(\lambda, \mu)$ and choosing the initial agents in $S_b(\lambda, \mu)$, thereby successfully constructing the state transition chain from $S_b(\lambda, \mu)$ to $S_g$. For further technical details, please refer to Lemma 5 and its proof.

**4.5 Step 3: Leftward Inertial Motion**

This section aims to construct a state transition chain from $S_b$ to $S_g$, corresponding to a single leftward inertia motion. On this basis, an upper bound $t_{bg}$ for the first hitting time from $S_b$ to $S_g$ and a lower bound $p_{bg}$ for the corresponding transition probability can be obtained. It is important to note that, owing to symmetry, the methods for constructing the state transition chain closely mirror those described in Section 4.3.

The state transition chains, along with the results of $t_{bg}$ and $p_{bg}$, are presented in Lemma 5, which is rigorously proven in Appendix A.11-Appendix A.13.

*Lemma 5*

When $\omega = 1$ and $pb \neq gb$, the following propositions hold:

(a) For the case where $ub - u_g > \frac{1}{3}d_0$ and $u_g \geq \frac{3}{4}pb + \frac{1}{4}gb$, given an initial state vector $\xi(t_0) \in S_b\left(\frac{1}{40}, \frac{1}{20}\right)$, we set $t_1 = \left\lceil \frac{ub - u_g}{|v(t_0)|} \right\rceil - 2$. The transition chain $\{\xi(t_0)\} \to S_1^3 \to S_2^3 \to \cdots S_{t_1-1}^3 \to S_g$ is feasible from $\{\xi(t_0)\}$ to $S_g$, where:

$$\begin{cases} S_t^3 = R_t^3 \times \mathbb{R}, t = 1, \ldots, t_1 - 1 \\ R_t^3 = [a_t - \delta_a, a_t]_I \\ a_t = ub - (t+1)|v(t_0)| - t(t-1)\delta_a \\ \delta_a = \frac{ub - u_g - (t_1+1)|v(t_0)|}{t_1(t_1-1)} \end{cases} \quad (4.27)$$

(b) For the case where $ub - u_g > \frac{1}{3}d_0$ and $u_g < \frac{3}{4}pb + \frac{1}{4}gb$, given an initial state vector $\xi(t_0) \in S_b\left(\frac{1}{40}, \frac{1}{20}\right)$, we set $t_2 = \left\lceil \frac{ub - \left(\frac{1}{4}pb + \frac{3}{4}gb\right)}{|v(t_0)|} \right\rceil - 2$. The transition chain $\{\xi(t_0)\} \to U_1^3 \to U_2^3 \to \cdots U_{t_2-1}^3 \to U_{t_2}^3$ is feasible from $\{\xi(t_0)\}$ to $U_{t_2}^3$, where:

$$\begin{cases} U_t^3 = W_t^3 \times \mathbb{R}, t = 1, \ldots, t_2 \\ W_t^3 = [b_t - \delta_b, b_t]_I \\ b_t = ub - (t+1)|v(t_0)| - t(t-1)\delta_b \\ \delta_b = \frac{ub - \frac{1}{4}pb - \frac{3}{4}gb - (t_2+1)|v(t_0)|}{t_2(t_2-1)} \end{cases} \quad (4.28)$$

Afterwards, for the state vector $\xi(t_0 + t_2) \in U_{t_2}^3$ within the transition chain $\{\xi(t_0)\} \to \cdots \to U_{t_2}^3$, we set $t_3 = \left\lceil \frac{x(t_0+t_2) - u_g}{|v(t_0+t_2)|} \right\rceil + 1$. The transition chain $\{\xi(t_0+t_2)\} \to X_1^3 \to X_2^3 \to \cdots X_{t_3-1}^3 \to S_g$ is feasible from $\{\xi(t_0+t_2)\}$ to $S_g$, where

$$\begin{cases} X_t^3 = Y_t^3 \times \mathbb{R}, t = 1, \ldots, t_3 - 1 \\ Y_t^3 = [c_t, c_t + \delta_c]_I \\ c_t = x(t_0 + t_2) - t|v(t_0 + t_2)| + t(t-1)\delta_c \\ \delta_c = \frac{u_g - x(t_0+t_2) + t_3|v(t_0+t_2)|}{t_3(t_3-1)} \end{cases} \quad (4.29)$$

(c) For the case $ub - u_g \leq \frac{1}{3}d_0$, it holds that $S_b\left(\frac{ub - u_g}{d_0}, \min\left\{\frac{ub - l_g}{d_0}, \frac{1}{2}\right\}\right) \subset S_g$; for the case where $ub - u_g > \frac{1}{3}d_0$, for any given state vector $\xi(t_0) \in S_b\left(\frac{1}{48}, \frac{1}{24}\right)$ at any initial iteration $t_0$, it holds that

$$\mathbb{P}\left\{T\left(S_g | \xi(t_0)\right) \leq t_{bg}\right\} \geq p_{bg} \quad (4.30)$$

where

$$\begin{cases} t_{bg} = \left\lceil 80\frac{ub-lb}{d_0} \right\rceil \\ p_{bg} = h^{t_{bg}}\left(10^{-5}\frac{d_0^3}{(ub-lb)^2}\right)h(|R_g|) \end{cases} \quad (4.31)$$

**4.6 Final proof of Proposition 2**

Sections 4.3 to 4.5 present the complete three-step construction of the state transition chain from $S_0$ to $S_g$. Based on the corresponding first hitting time and transition probability results from Lemmas 2, 4, and 5, we can provide explicit expressions for $t_{e0}$ and $p_{e0}$ along with Eq. (3.1), thus completing the proof of Proposition 2 (see Appendix A.14 for the proof details).

*Proposition 2 (with explicit expressions for $t_{e0}$ and $p_{e0}$)*

When $\omega = 1$ and $pb \neq gb$, there exist $p_{e0} > 0$ and $t_{e0} \in \mathbb{N}^+$ such that for any given state vector $\xi(t_0)$ at any initial iteration $t_0$, it holds that

$$\mathbb{P}\left\{T\left(S_g|\xi(t_0)\right) \leq t_{e0}\right\} \geq p_{e0} \quad (4.32)$$

where

$$\begin{cases} t_{e0} = \left\lceil \frac{2(C_1+C_2)+100(ub-lb)}{d_0} \right\rceil \\ p_{e0} = h^{t_{e0}}\left(10^{-5}\frac{d_0^3}{(ub-lb)^2}\right)h(|R_g|) \\ d_0 = \min\{C_1, C_2, 1\}|gb - pb| \end{cases} \quad (4.33)$$

where $R_g$ represents the neighborhood of the GO, and $d_0$ is defined in Eq. (4.4).

Finally, we summarize the process of proving Proposition 2 described above in Fig. 6.

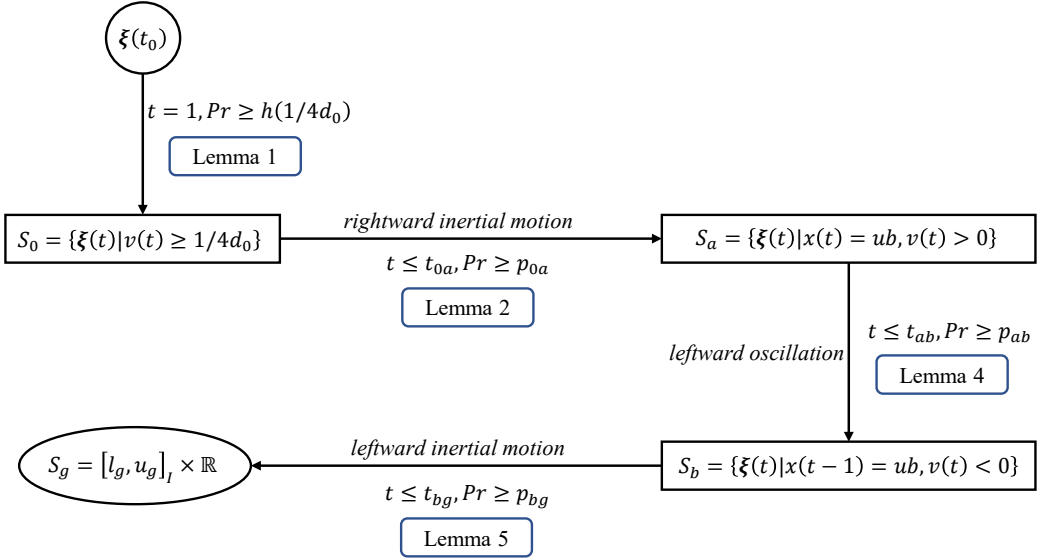

Fig. 6 Proof process of Proposition 2

**4.7 Discussions**

Upon reviewing the preceding analysis and the values of $t_{e0}$ and $p_{e0}$ in Eq. (4.33), the following observations can be made:

(1) If $pb = gb$, then $d_0 = 0$, resulting in $t_{e0} = +\infty$ and $p_{e0} = 0$. This indicates that the analysis is meaningful only when $pb \neq gb$, meaning that the agent's pbest and the population's gbest stagnate near different local optima. According to Corollary 1, $pb \neq gb$ ensures that $d_0 > 0$, thus guaranteeing that $v(t+1)$ has a certain search range regardless of a fixed $\xi(t)$, and consequently, assuring the ability of PSO to escape from LO. If $pb = gb$ and $v(t_0) = 0$, it follows that $v(t_0 + 1) = v(t_0 + 2) = \cdots = 0$, indicating that agent's position remains unchanged, resulting in a loss of its search capability within the feasible region. Therefore, when the pbest and gbest stagnate near different local optima, the agent retains the ability to escape from LO. Conversely, when they stagnate near the same LO, the agent may lose this capability. In addition, as noted from (40), as the value of $gb - pb$ increases, $t_{e0}$ decreases and $p_{e0}$ increases, suggesting that a greater distance between pbest and gbest enhances the ability of agent to escape from LO.

(2) The inertial weight $\omega = 1$ is a crucial condition for ensuring the validity of the above analysis. For instance, in the construction of the rightward inertial motion addressed in Section 5.3, it is precisely because of $\omega = 1$ that the state transition chain can be constructed near the arithmetic sequence with an initial term $x(0)$ and a common difference $v(0)$. If $\omega < 1$, as illustrated in Case 5, the agent may fail to reach the upper bound $ub$ within a single rightward inertial motion, indicating a reduction in the agent's search ability and a lack of guaranteed ability to escape from LO.

(3) As observed from Eq. (4.33), increasing the values of learning factors $C_1$ and $C_2$ contributes to an increase in $d_0$, which in turn facilitates a decrease in $t_{e0}$ and an increase in $p_{e0}$. This suggests that increasing the values of $C_1$ and $C_2$ enhances the ability of PSO to escape from LO.

## 5. NUMERICAL EXPERIMENTS

In this section, several numerical experiments are conducted to validate the conclusions presented above.

First, we perform experiments to assess the ability of PSO to escape from LO under the condition $\omega = 1$. We set $lb = 0, ub = 20, C_1, C_2 \in \{2, 2.5\}, x(0) \sim U(0,2), v(0) \sim U(-1,1)$, with the stagnated $pb = 3$ and $gb \in \{3,4,5\}$, while the neighborhood of the GO is set as $R_g = [19, 20]_I$. Since $x(0), pb$ and $gb$ are distant from $R_g$, combined with a small initial agent velocity $|v(0)|$, thus finding a position within $R_g$ becomes difficult, and the ability of PSO to escape from LO can be adequately validated through this experiment. For each set of $(C_1, C_2, gb)$, PSO is independently executed $10^5$ times, and the proportion of successful runs that escape from the LO within $t$ iterations allows us to approximate $\mathbb{P}\{T(S_g|\xi(0)) \leq t\}$, as illustrated in Fig. 7.

From Fig. 7, as $t$ increases, the value of $\mathbb{P}\{T(S_g|\xi(0)) \leq t\}$ converges to 1 if and only if $pb \neq gb$, thereby confirming Proposition 3. Additionally, it is observed that $\mathbb{P}\{T(S_g|\xi(0)) \leq t\}$ increases with $|pb - gb|, C_1$ and $C_2$, indicating that a larger the distance between the stagnated pbest and gbest, along with greater values of $C_1$ and $C_2$, enhances PSO's ability to escape from LO.

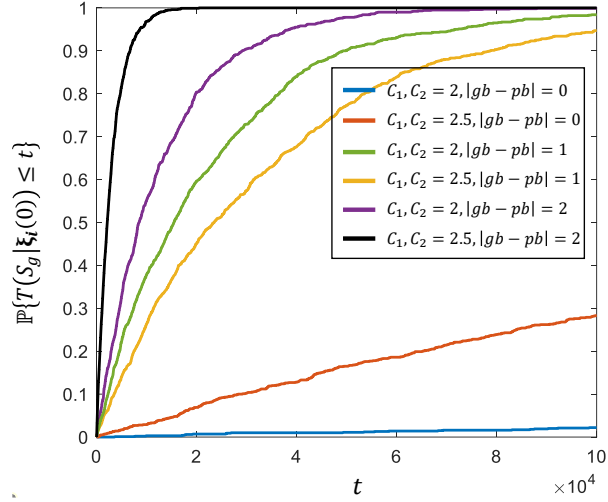

Fig. 7 Approximate values of $\mathbb{P}\{T(S_g|\boldsymbol{\xi}(0)) \le t\}$.

Next, we conduct the experiments to assess the ability of PSO to escape from LO under the condition $\omega < 1$. Here, we set $lb = 0, ub \in \{20,22,24,26,28,30\}$ and $R_g = [ub-1, ub]_I$. In PSO, we specify $C_1, C_2 \in \{1.6, 2, 2.4\}$, $\omega \in \{0.7, 0.8, 0.9\}$, $x(0) \sim U(0,2), v(0) \sim U(-1,1)$, $pb = 3$ and $gb = 4$. For each set of $(C_1, C_2, \omega, ub)$, PSO is executed independently $10^4$ times. By calculating the ratio of successful escapes from the LO over $10^7$ iterations in these $10^4$ runs, we approximate the value of $p_e$, with the corresponding experimental results presented in Table I. From Table I, it is evident that $p_e = 0$ holds for sufficiently large distance between $x(0)$ and $R_g$, and for a given $R_g$, the value of $p_e$ increases with $C_1, C_2$, and $\omega$, effectively validating the aforementioned conclusions.

TABLE I Approximate values of $p_e$

| $\omega$ | $C_1, C_2$ | $R_g$ | | | | | |
|---|---|---|---|---|---|---|---|
| | | $[19,20]_I$ | $[21,22]_I$ | $[23,24]_I$ | $[25,26]_I$ | $[27,28]_I$ | $[29,30]_I$ |
| 0.9 | 2.4 | 1 | 1 | 1 | 1 | 1 | 0.71 |
| | 2 | 1 | 1 | 1 | 0.63 | 0.25 | 0.06 |
| | 1.6 | 1 | 0.9 | 0.36 | 0.06 | 0 | 0 |
| 0.8 | 2.4 | 1 | 1 | 1 | 1 | 0.77 | 0.10 |
| | 2 | 1 | 1 | 0.48 | 0.10 | 0.02 | 0 |
| | 1.6 | 0.83 | 0.09 | 0.01 | 0.01 | 0 | 0 |
| 0.7 | 2.4 | 1 | 1 | 1 | 0.85 | 0.06 | 0.01 |
| | 2 | 1 | 0.62 | 0.02 | 0 | 0 | 0 |
| | 1.6 | 0.03 | 0.01 | 0 | 0 | 0 | 0 |

Finally, the Proposition 3 will be interpreted from the perspective of probability distribution. It was discussed in Section 3 that, as Proposition 3 implies, the search range of the agent can cover the entire feasible region after several iterations, ensuring that the probability of escaping from the LO in each iteration is greater than zero. By setting $[lb, ub]_I = [0,9], pb = 0.5, gb = 2, \boldsymbol{\xi}(0) = [1, 0.5], R_g = [8.5, 9]_I, \omega = 1$, and $C_1 = C_2 = 2$, we apply Monte Carlo simulations to derive the distribution of $x(t)$, as shown in Fig. 8. The distribution of $x(t)$ within the open interval $(lb, ub)$ is depicted in Fig. 8(a), while the simulated values of $\mathbb{P}\{x(t) = lb\}, \mathbb{P}\{x(t) = ub\}$, and $\mathbb{P}\{x(t) \in R_g\}$ are depicted in Fig. 8(b). As shown in Fig. 8, although $x(0)$ is far from $R_g$ and the initial velocity $v(0)$ is small, the distribution range of $x(t)$ continuously expands through iterations, ultimately covering $[lb, ub]_I$ (and therefore $R_g$) after the fifth iteration. Subsequently, the distribution of $x(t)$ and the value of $\mathbb{P}\{x(t) \in R_g\}$ tend to stabilize, thereby ensuring a positive probability of escaping from the LO in each iteration.

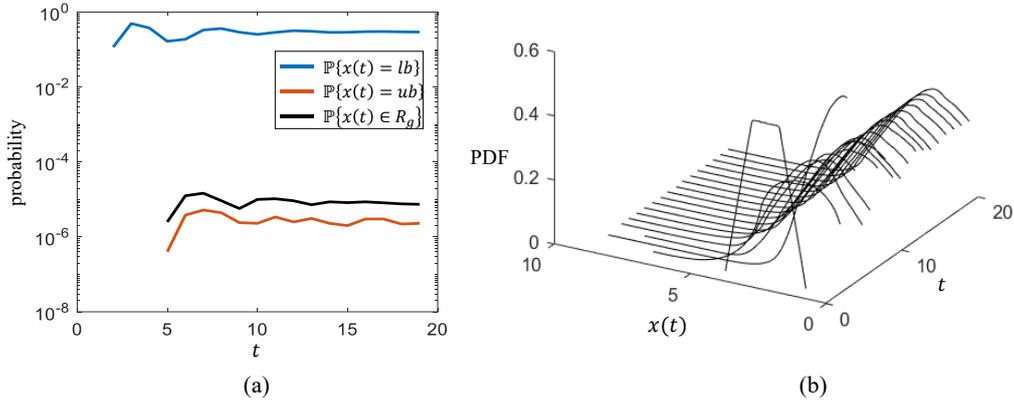

Fig. 8. Simulation results of the distribution of $x(t)$. (a) Distribution of $x(t)$ within $(lb, ub)$. (b) Simulated values of $\mathbb{P}\{x(t) = lb\}, \mathbb{P}\{x(t) = ub\}$ and $\mathbb{P}\{x(t) \in R_g\}$.

# 6. CONCLUSION

This work proposes a theoretical analysis framework to investigate the ability and mechanism of PSO in escaping from LO. By jointly modelling the position and velocity of each search agent as a continuous state Markov chain, we rigorously calculate the probability that PSO escapes from LO within a finite number of iterations and carefully analyze the behavior of agents prior to this escape. The mathematical calculations reveal that, for ensuring the probability of escaping from LO within finite iterations to be 1, the conditions that pbest and gbest stagnate around different local optima and the inertial weight $\omega = 1$ in PSO are both necessary and sufficient. In addition, the theoretical analyses also reveal that increasing the distance between pbest and gbest, as well as the values of $C_1$ and $C_2$, can enhance the speed at which PSO escapes from LO. Finally, several numerical experiments are conducted to validate the theoretical conclusions presented. In summary, this work addresses a gap in the theoretical analysis of PSO's ability to escape from LO, providing a more detailed and fundamental characterization of PSO's global search capacity and mechanisms, as well as offering a solid theoretical basis for the development of effective improvement strategies for PSO.

# APPENDIX

### A.1 Proof of Proposition 3 (based on Proposition 2)

For $n \in \mathbb{N}^+$, according to Proposition 2, it holds that

$$\mathbb{P}\{T(S_g|\xi(t_0)) > (n+1)t_{e0} | T(S_g|\xi(t_0)) > nt_{e0}\} \leq 1 - p_{e0}$$

Therefore,

$$\mathbb{P}\{T(S_g|\xi(t_0)) > (n+1)t_{e0}\}$$

$$= \mathbb{P}\{T(S_g|\xi(t_0)) > (n+1)t_{e0} | T(S_g|\xi(t_0)) > nt_{e0}\} \cdot \mathbb{P}\{T(S_g|\xi(t_0)) > nt_{e0}\}$$

$$\leq (1 - p_{e0})\mathbb{P}\{T(S_g|\xi(t_0)) > nt_{e0}\}$$

This implies that

$$\mathbb{P}\{T(S_g|\xi(t_0)) > nt_{e0}\} \leq (1 - p_{e0})^n$$

As a result,

$$p_e = \lim_{n\to\infty}\left[1 - \mathbb{P}\{T(S_g|\xi(t_0)) > nt_e\}\right] \geq 1 - \lim_{n\to\infty}(1 - p_{e0})^n = 1$$

□

### A.2 Supplement for Proposition 4

The values of $v_{fi}$ and $h_f$ are as follows:

$$[v_{fi}]_{i=1,\dots,4} = v(t) + \begin{cases} [0, dv_1, dv_2, dv_1 + dv_2], x(t) \leq pb \\ [dv_1 + dv_2, dv_1, dv_2, 0], x(t) \geq gb \\ [dv_1, 0, dv_1 + dv_2, dv_2], pb < x(t) < gb, dv_1 + dv_2 > 0 \\ [dv_1, dv_1 + dv_2, 0, dv_2], pb < x(t) < gb, dv_1 + dv_2 \leq 0 \end{cases} \quad (A1)$$

$$h_f = \frac{2}{v_{f4} + v_{f3} - v_{f2} - v_{f1}} \quad (A2)$$

where

$$\begin{cases} dv_1 = \min\{C_1(pb - x(t)), C_2(gb - x(t))\} \\ dv_2 = \max\{C_1(pb - x(t)), C_2(gb - x(t))\} \end{cases} \quad (A3)$$

### A.3 Proof of Proposition 4

It has been pointed out that $v(t+1) = [v(t) + C_1r_1(pb - x(t))] + C_2r_2(gb - x(t))$ can be expressed as the sum of two random variables following uniform distributions. Here, we set $X \sim U(l_x, u_x), Y \sim U(l_y, u_y)$, then, by calculating the PDF of $Z = X + Y$, the PDF of $v(t+1)$ can be easily obtained.

By introducing the unit step function $\varepsilon$, the PDFs of $X$ and $Y$ can be expressed as

$$f_X(x) = \frac{1}{u_x-l_x}[\varepsilon(x-l_x) - \varepsilon(x-u_x)], f_Y(y) = \frac{1}{u_y-l_y}[\varepsilon(y-l_y) - \varepsilon(y-u_y)]$$

On this basis, according to the property of convolution operator, the PDF of $Z$ can be expressed as $f_Z(z) = \int_{\mathbb{R}} f_X(x) f_Y(z-x) dx$.

For the case $l_x + l_y \leq z \leq \min\{l_x + u_y, u_x + l_y\}$, it can be derived that:

$$f_Z(z) = \int_{l_x}^{z-l_y} f_X(x) f_Y(z-x) dx = \frac{1}{(u_x-l_x)(u_y-l_y)} [z - (l_x + l_y)]$$

For the case $\min\{l_x + u_y, u_x + l_y\} \leq z \leq \max\{l_x + u_y, u_x + l_y\}$, it can be derived that:

$$f_Z(z) = \begin{cases} \int_{l_x}^{u_x} f_X(x) f_Y(z-x) dx, u_y - l_y \geq u_x - l_x \\ \int_{z-u_y}^{z-l_y} f_X(x) f_Y(z-x) dx, u_y - l_y \leq u_x - l_x \end{cases}$$

$$= \frac{1}{(u_x-l_x)(u_y-l_y)} [\min\{l_x + u_y, u_x + l_y\} - (l_x + l_y)]$$

For the case $\max\{l_x + u_y, u_x + l_y\} \leq z \leq u_x + u_y$, it can be derived that

$$f_Z(z) = \int_{z-u_y}^{u_x} f_X(x) f_Y(z-x) dx = \frac{1}{(u_x-l_x)(u_y-l_y)} [-z + (u_x + u_y)].$$

To sum up, it holds that

$$f_Z(z) = \begin{cases} \frac{z-(l_x+l_y)}{(u_x-l_x)(u_y-l_y)}, l_x + l_y \leq z \leq \min\{l_x + u_y, u_x + l_y\} \\ \frac{\min\{l_x+u_y,u_x+l_y\}-(l_x+l_y)}{(u_x-l_x)(u_y-l_y)}, \min\{l_x + u_y, u_x + l_y\} \leq z \leq \max\{l_x + u_y, u_x + l_y\} \\ \frac{-z+(u_x+u_y)}{(u_x-l_x)(u_y-l_y)}, \max\{l_x + u_y, u_x + l_y\} \leq z \leq u_x + u_y \\ 0, otherwise \end{cases}$$

Finally, by substituting the upper bounds and lower bounds of random variables $v(t) + C_1 r_1(pb - x(t))$ and $C_2 r_2(gb - x(t))$ into the expression of $f_Z(z)$, Eq. (4.3) is obtained.

□

### A.4 Proof of Corollary 1

*Proof of Corollary 1(a)*

According to Proposition 4, $D_{v(t+1)} = [v_{f1}, v_{f4}]$ and $v(t) \in D_{v(t+1)}$ obviously hold.

□

*Proof of Corollary 1(b)*

According to Eqs. (A1)-(A3), for the case $x(t) < pb$, it holds that:
$|D_{v(t+1)}| = C_1(pb - x(t)) + C_2(gb - x(t)) \geq C_2(gb - pb) \geq \min\{C_1, C_2\}(gb - pb)$
For the case $x(t) > gb$, it holds that:
$|D_{v(t+1)}| = C_1(x(t) - pb) + C_2(x(t) - gb) \geq C_1(gb - pb) \geq \min\{C_1, C_2\}(gb - pb)$
For the case $pb \leq x(t) \leq gb$, it holds that:
$|D_{v(t+1)}| = C_2(gb - x(t)) - C_1(pb - x(t)) \geq \min\{C_2(gb - pb), -C_1(pb - gb)\} = \min\{C_1, C_2\}(gb - pb)$

To sum up, it holds that $|D_{v(t+1)}| \geq \min\{C_1, C_2\}(gb - pb) \geq \min\{C_1, C_2, 1\}(gb - pb)$.

□

*Proof of Corollary 1(c)*

First, we estimate the lower bound of $\mathbb{P}\{v(t+1) \in D_1\}$ by providing a lower bound of $f_{v(t+1)}(v|\xi(t))$ for $v \in D_{v(t+1)}$. Specifically, if $\forall v \in D_{v(t+1)}, g(v) \leq f_{v(t+1)}(v|\xi(t))$, then it holds that $\mathbb{P}\{v(t+1) \in D_1\} = \int_{v \in D_1} f_{v(t+1)}(v|\xi(t))dv \geq \int_{v \in D_1} g(v)dv$.

For the cases $D_1 \subset [v_{f1}, v_{f2}]$ or $D_1 \subset [v_{f3}, v_{f4}]$, $g(v)$ is shown as the blue curve in Figure A1(a) (where the black curve denotes the curve of $f_{v(t+1)}(v|\xi(t))$). On this basis, it holds that:

$$\mathbb{P}\{v(t+1) \in D_1\} \geq \int_{v \in D_1} g(v)dv = \frac{1}{2}\frac{h_f}{v_{f2}-v_{f1}}|D_1|^2 \tag{A4}$$

For the case $D_1 \cap [v_{f2}, v_{f3}] \neq \Phi$, $g(v)$ is shown as the blue curve in Figure A1(b) (where the black curve denotes the curve of $f_{v(t+1)}(v|\xi(t))$). On this basis, it holds that:

$$\mathbb{P}\{v(t+1) \in D_1\} \geq \int_{v \in D_1} g(v)dv = \frac{1}{2}h_f|D_1| \tag{A5}$$

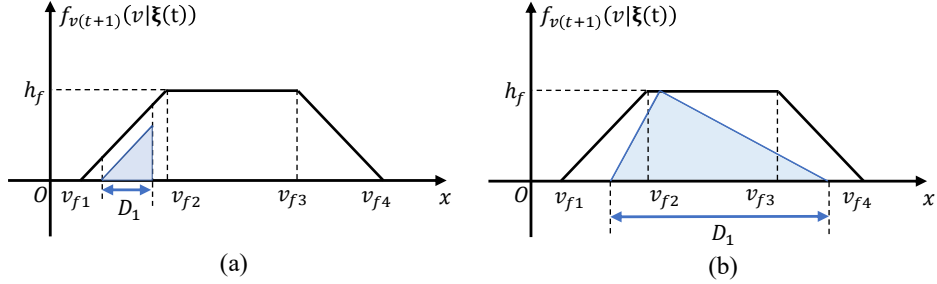

Figure A1: curves of $f_{v(t+1)}(v|\xi(t))$ and $g(v)$. (a): $D_1 \subset [v_{f1}, v_{f2}]$ or $D_1 \subset [v_{f3}, v_{f4}]$, (b) $D_1 \cap [v_{f2}, v_{f3}] \neq \Phi$.

Next, we find the lower bounds for $\frac{h_f}{v_{f2}-v_{f1}}$ and $h_f$. According to Eqs. (A1)-(A3), for the case $x(t) \leq pb$, it holds that:

$$h_f = \frac{2}{v_{f4}+v_{f3}-v_{f2}-v_{f1}} = \frac{1}{\max\{C_1(pb-x(t)), C_2(gb-x(t))\}} \geq \frac{1}{\max\{C_1, C_2\}(ub-lb)}$$

$$\frac{h_f}{v_{f2}-v_{f1}} = \frac{1}{C_1 C_2 (pb-x(t))(gb-x(t))} \geq \frac{1}{C_1 C_2 (ub-lb)^2}$$

For the case $x(t) \geq gb$, it holds that:

$$h_f = \frac{2}{v_{f4}+v_{f3}-v_{f2}-v_{f1}} = \frac{1}{\max\{C_1(x(t)-pb), C_2(x(t)-gb)\}} \geq \frac{1}{\max\{C_1, C_2\}(ub-lb)}$$

$$\frac{h_f}{v_{f2}-v_{f1}} = \frac{1}{C_1 C_2 (x(t)-pb)(x(t)-gb)} \geq \frac{1}{C_1 C_2 (ub-lb)^2}$$

For the case $pb \leq x(t) \leq gb$, it holds that:

$$h_f = \frac{2}{v_{f4}+v_{f3}-v_{f2}-v_{f1}} = \frac{1}{\max\{C_1(x(t)-pb), C_2(gb-x(t))\}} \geq \frac{1}{\max\{C_1, C_2\}(ub-lb)}$$

$$\frac{h_f}{v_{f2}-v_{f1}} = \frac{1}{C_1 C_2 (x(t)-p)(g-x(t))} \geq \frac{1}{\frac{1}{4}C_1 C_2 (g-p)^2} \geq \frac{4}{C_1 C_2 (ub-lb)^2}.$$

To sum up, it holds that

$$\begin{cases} \inf\limits_{pb,gb,x(t),v(t)} \frac{1}{2}h_f \geq \frac{1}{2\max\{C_1,C_2\}(ub-lb)} \\ \inf\limits_{pb,gb,x(t),v(t)} \frac{1}{2}\frac{h_f}{v_{f2}-v_{f1}} \geq \frac{1}{2C_1C_2(ub-lb)^2} \end{cases} \quad (A6)$$

Finally, combining Eqs. (A4) to (A6) yields Eqs. (4.6) to (4.7).

□

## A.5 Proof of Lemma 1

Since $\max\{-v_{f1}, v_{f4}\} \geq \frac{1}{2}d_0$ (otherwise $|D_{v(t+1)}| = v_{f4} - v_{f1} < d_0$ contradicts with Eq. (4.5)), we only need to consider the cases $v_{f4} \geq \frac{1}{2}d_0$ and $v_{f1} \leq -\frac{1}{2}d_0$.

For the case $v_{f1} \leq -\frac{1}{2}d_0$, according to Corollary 1(c), it holds that

$$\mathbb{P}\left\{v(t+1) \geq \frac{1}{4}d_0 \text{ or } v(t+1) \leq -\frac{1}{4}d_0\right\}$$

$$\geq \mathbb{P}\left\{v(t+1) \leq -\frac{1}{4}d_0\right\}$$

$$\geq \mathbb{P}\left\{v_{f1} \leq v(t+1) \leq v_{f1} + \frac{1}{4}d_0\right\} \geq h\left(\frac{1}{4}d_0\right)$$

For the case $v_{f4} \geq \frac{1}{2}d_0$, according to Corollary 1(c), it holds that

$$\mathbb{P}\left\{v(t+1) \geq \frac{1}{4}d_0 \text{ or } v(t+1) \leq -\frac{1}{4}d_0\right\}$$

$$\geq \mathbb{P}\left\{v(t+1) \geq \frac{1}{4}d_0\right\}$$

$$\geq \mathbb{P}\left\{v_{f4} - \frac{1}{4}d_0 \leq v(t+1) \leq v_{f4}\right\} \geq h\left(\frac{1}{4}d_0\right)$$

To sum up, Eq. (4.10) obviously holds.

□

## A.6 Proof of Lemma 2(a)

We begin by proving the following Lemma A1:

*Lemma A1*

Given $x(t-1)$ and $x(t)$, and setting the distribution range of $x(t+1)$ as $D_{x(t+1)}$, it follows that, when neglecting boundary constraints, $D_{x(t+1)} \supset E_{x(t+1)}$, where

$$E_{x(t+1)} = 2x(t) - x(t-1) + \begin{cases} [-C_1(1-\lambda)(gb-pb), 0]_I, x(t) > \lambda pb + (1-\lambda)gb \\ [0, C_2\lambda(gb-pb)]_I, x(t) \leq \lambda pb + (1-\lambda)gb \end{cases} \quad (A7)$$

*Proof of Lemma A1*

When neglecting boundary constraints, substituting Eq. (2.8) into Eq. (2.7) yields:

$$x(t+1) = 2x(t) - x(t-1) + C_1r_1(pb - x(t)) + C_2r_2(gb - x(t)) \quad (A8)$$

When $x(t) > \lambda pb + (1-\lambda)gb$, setting $r_1 = r_2 = 0$ in Eq. (A8) yields $x(t+1) = 2x(t) - x(t-1)$, setting $r_1 = 1, r_2 = 0$ in Eq. (A8) yields $x(t+1) = 2x(t) - x(t-1) +$ $C_1(pb - x(t)) < 2x(t) - x(t-1) - C_1(1-\lambda)(gb - pb)$. Since $D_{x(t+1)}$ is a close interval (see Corollary 1(a)), it holds that:

$$D_{x(t+1)} \supset [2x(t) - x(t-1) + C_1(pb - x(t)), 2x(t) - x(t-1)]$$
$$\supset [2x(t) - x(t-1) - C_1(1-\lambda)(gb - pb), 2x(t) - x(t-1)] = E_{x(t+1)}$$

When $x(t) \leq \lambda pb + (1-\lambda)gb$, setting $r_1 = r_2 = 0$ in Eq. (A8) yields $x(t+1) = 2x(t) - x(t-1)$, and setting $r_1 = 0, r_2 = 1$ in Eq. (A8) yields $x(t+1) = 2x(t) - x(t-1) + C_2(gb - x(t)) \geq 2x(t) - x(t-1) + C_2\lambda(gb - pb)$. Since $D_{x(t+1)}$ is a close interval (see Corollary 1(a)), it holds that:

$$D_{x(t+1)} \supset [2x(t) - x(t-1), 2x(t) - x(t-1) + C_2(gb - x(t))]$$
$$\supset [2x(t) - x(t-1), 2x(t) - x(t-1) + C_2\lambda(gb - pb)] = E_{x(t+1)}$$

□

Returning to the proof of Lemma 2(a). We set $t_1$ as the minimum $t$ that satisfies $\inf R_t^1 \geq ub$, then it holds that

$$\begin{cases} a_{t_1} - \delta_a \geq ub \\ a_{t_1-1} - \delta_a < ub \end{cases} \tag{A9}$$

On this basis, to prove the feasibility of the state transition chain $\{\xi(t_0)\} \to S_1^1 \to S_2^1 \to \cdots S_{t_1-1}^1 \to S_a$ for $\xi(t_0) \in S_0$ satisfying $x(t_0) \geq \lambda pb + (1-\lambda)gb$, we will sequentially introduce and prove the following three propositions.

*Proposition A6.1*

$R_t^1\left(t = 0, \ldots, \left\lceil \frac{8(ub-lb)}{d_0} \right\rceil - 1\right)$ represents a series of increasing intervals. Specifically, it holds that $\forall t = 1, \ldots, \left\lceil \frac{8(ub-lb)}{d_0} \right\rceil, \inf R_t^1 > \sup R_{t-1}^1$.

*Proof of Proposition A6.1*

For any $t = 0, \ldots, \left\lceil \frac{8(ub-lb)}{d_0} \right\rceil - 1$, it holds that

$\inf R_t^1 - \sup R_{t-1}^1$
$= a_t - \delta_a - a_{t-1}$
$= v(t_0) - (2t-1)\delta_a$
$> \frac{1}{4}d_0 - \left(2\frac{8(ub-lb)}{d_0} - 1\right)\frac{d_0^2}{128(ub-lb)}$
$= \frac{1}{8}d_0 + \frac{d_0^2}{128(ub-lb)} > \frac{1}{8}d_0$

□

*Proposition A6.2*

The $t_1$ defined by Eq. (A9) exists and is unique.

*Proof of Proposition A6.2*

According to Proposition A6.1, $R_t^1\left(t = 1, \ldots, \left\lceil \frac{8(ub-lb)}{d_0} \right\rceil - 1\right)$ represents a series of increasing intervals. Thus, we only need to prove that $a_{\left\lceil \frac{8(ub-lb)}{d_0} \right\rceil - 1} - \delta_a \geq ub$. According to Proposition A6.1, it holds that:

$$a_{\left\lceil\frac{8(ub-lb)}{d_0}\right\rceil-1} - \delta_a$$

$$= a_1 + \sum_{t=2}^{\left\lceil\frac{8(ub-lb)}{d_0}\right\rceil-1}(a_t - a_{t-1}) - \delta_a$$

$$\geq a_1 + \sum_{t=2}^{\left\lceil\frac{8(ub-lb)}{d_0}\right\rceil-1}(\inf R_t^1 - \sup R_{t-1}^1 + \delta_a) - \delta_a$$

$$\geq x(t_0) + v(t_0) + \sum_{t=2}^{\left\lceil\frac{8(ub-lb)}{d_0}\right\rceil-1}\left(\frac{1}{8}d_0\right)$$

$$\geq x(t_0) + \frac{1}{4}d_0 + \frac{1}{8}d_0\left(\left\lceil\frac{8(ub-lb)}{d_0}\right\rceil - 2\right)$$

$$= x(t_0) + \frac{1}{8}d_0\left\lceil\frac{8(ub-lb)}{d_0}\right\rceil$$

$$> x(t_0) + ub - lb \geq ub$$

□

*Proposition A6.3*

For any $t = 0, \ldots, t_1 - 1$, any $x(t + t_0 - 1) \in R_{t-1}^1, x(t + t_0) \in R_t^1$, the distribution range of $x(t + t_0 + 1)$ covers $R_{t+1}^1$.

*Proof of Proposition A6.3*

According to Lemma A1, we only need to prove that $\forall t = 0, \ldots, t_1 - 1, \forall x(t + t_0 - 1) \in R_{t-1}^1, \forall x(t + t_0) \in R_t^1, E_{x(t+t_0+1)} \supset R_{t+1}^1$ holds. According to Proposition A6.1 and Eq. (A9), $\forall x(t + t_0) \in R_t^1$, it holds that $x(t + t_0) \geq \inf R_0^1 \geq \lambda pb + (1 - \lambda)gb$, thus the expression of $E_{x(t+t_0+1)}$ becomes (see Eq. (A7))

$$E_{x(t+t_0+1)} = [2x(t+t_0) - x(t+t_0-1) - C_1(1-\lambda)(gb-pb), 2x(t+t_0) - x(t+t_0-1)]_I \tag{A10}$$

Therefore, we only need to prove that

$$\begin{cases} \inf_{x(t+t_0-1)\in R_{t-1}^1, x(t+t_0)\in R_t^1}[2x(t+t_0) - x(t+t_0-1)] \geq a_{t+1} \\ \sup_{x(t+t_0-1)\in R_{t-1}^1, x(t+t_0)\in R_t^1}[2x(t+t_0) - x(t+t_0-1) - C_1(1-\lambda)(gb-pb)] \leq a_{t+1} - \delta_a \end{cases} \tag{A11}$$

Here,

$$\inf_{x(t+t_0-1)\in R_{t-1}^1, x(t+t_0)\in R_t^1}[2x(t+t_0) - x(t+t_0-1)] - a_{t+1}$$

$$\geq 2(a_t - \delta_a) - a_{t-1} - a_{t+1} = 0$$

and

$$\sup_{x(t+t_0-1)\in R_{t-1}^1, x(t+t_0)\in R_t^1}[2x(t+t_0) - x(t+t_0-1) - C_1(1-\lambda)(gb-pb)] - (a_{t+1} - \delta_a)$$

$$\leq 2a_t - (a_{t-1} - \delta_a) - C_1(1-\lambda)(gb-pb) - (a_{t+1} - \delta_a)$$
$$= -C_1(1-\lambda)(gb-pb) + 4\delta_a$$
$$\leq -(1-\lambda)d_0 + 4\frac{1-\lambda}{4}d_0 = 0$$

□

By combining Propositions A6.1-A6.3, we successfully verify the feasibility of the state transition chain for agent's position $R_0^1 \to R_1^1 \to \cdots \to R_{t_1}^1$ when neglecting boundary constraints. Finally, we apply the boundary constraints to complete the proof of Lemma 2(a).

According to Eq. (A9), $t_1$ satisfies $\inf R_{t_1}^1 \geq ub$ and $\inf R_{t_1-1}^1 < ub$.

In the case where $\sup R_{t_1-1}^1 < ub$, the agent's position does not reach $ub$ unless it reaches the interval $R_{t_1}^1$. Thus, the feasibility of the state transition chain $\{\xi(t_0)\} \to S_1^1 \to S_2^1 \to \cdots S_{t_1-1}^1 \to S_a$ is verified, with first hitting time $t_{0a} = t_1$.

In the case where $\sup R_{t_1-1}^1 \geq ub$, when the agent's position $x(t)$ reaches the interval $R_{t_1-1}^1$, both $x(t) < ub$ and $x(t) \geq ub$ become possible. If $x(t) < ub$, the corresponding first hitting time of the state transition chain becomes $t_1$; if $x(t) \geq ub$, then it is constrained to $ub$, indicating that the agent reaches $S_a$ at iteration $t_1 - 1$, in which case the first hitting time becomes $t_1 - 1$.

In summary, we have completed the proof of Lemma 2(a).

□

### A.7 Proof of Lemma 2(b)

We further define $W_0^1 = \{x(t_0)\}$ and set $t_2$ as the minimum $t$ such that $\sup W_t^1 \geq \frac{1}{4}pb + \frac{3}{4}gb$. Thus, it holds that

$$\begin{cases} b_{t_2-1} + \delta_b < \frac{1}{4}pb + \frac{3}{4}gb \\ b_{t_2} + \delta_b \geq \frac{1}{4}pb + \frac{3}{4}gb \end{cases} \quad (A12)$$

On this basis, to prove the feasibility of the state transition chain $\{\xi(t_0)\} \to U_1^1 \to U_2^1 \to \cdots U_{t_2-1}^1 \to U_m$ for $\xi(t_0) \in S_0$ satisfying $x(t_0) < \frac{1}{4}pb + \frac{3}{4}gb$, we will sequentially introduce and prove the following three propositions.

*Proposition A7.1*

$W_t^1\left(t = 0, \ldots, \left\lceil \frac{4(ub-lb)}{d_0} \right\rceil\right)$ represents a series of increasing intervals. Specifically, it holds that $\forall t = 1, \ldots, \left\lceil \frac{4(ub-lb)}{d_0} \right\rceil, \inf W_t^1 > \sup W_{t-1}^1$.

*Proof of Proposition A7.1*

For $t = 1$, $\inf W_t^1 - \sup W_{t-1}^1 = v(t_0) > 0$; For $t = 2, \ldots, \left\lceil \frac{4(ub-lb)}{d_0} \right\rceil$, it holds that

$\inf W_t^1 - \sup W_{t-1}^1$
$= b_t - (b_{t-1} + \delta_b)$
$= v(t_0) + (2t - 3)\delta_b$
$> v(t_0) > 0$

□

*Proposition A7.2*

The $t_2$ defined by Eq. (A12) exists and is unique.

*Proof of Proposition A7.2*

According to Proposition A7.1, $W_t^1\left(t = 1, \ldots, \left\lceil \frac{4(ub-lb)}{d_0} \right\rceil\right)$ represents a series of increasing intervals. Thus, we only need to prove that $b_{\left\lceil \frac{4(ub-lb)}{d_0} \right\rceil} + \delta_b > \frac{1}{4}pb + \frac{3}{4}gb$. According to Proposition A7.1, it holds that:

$$b_{\left\lceil \frac{4(ub-lb)}{d_0} \right\rceil} + \delta_b$$

$$= b_1 + \sum_{t=2}^{\left\lceil \frac{4(ub-lb)}{d_0} \right\rceil}(b_t - b_{t-1}) + \delta_b$$

$$\geq b_1 + \sum_{t=2}^{\left\lceil \frac{4(ub-lb)}{d_0} \right\rceil}[v(t_0) + (2t - 2)\delta_b] + \delta_b$$

$$> x(t_0) + v(t_0) + \sum_{t=2}^{\left\lceil \frac{4(ub-lb)}{d_0} \right\rceil}\left(\frac{1}{4}d_0\right)$$

$$\geq x(t_0) + \frac{1}{4}d_0 + \frac{1}{4}d_0 \left(\left\lceil \frac{4(ub-lb)}{d_0} \right\rceil - 1\right)$$

$$\geq x(t_0) + ub - lb$$

$$\geq ub > \frac{1}{4}pb + \frac{3}{4}gb$$

□

*Proposition A7.3*

For any $t = 0, \ldots, t_2 - 1$, any $x(t + t_0 - 1) \in W_{t-1}^1, x(t + t_0) \in W_t^1$, the distribution range of $x(t + t_0 + 1)$ covers $W_{t+1}^1$.

*Proof of Proposition A7.3*

According to Lemma A1, we only need to prove that $\forall t = 0, \ldots, t_2 - 1, \forall x(t + t_0 - 1) \in W_{t-1}^1, \forall x(t + t_0) \in W_t^1, E_{x(t+t_0+1)} \supset W_{t+1}^1$ holds. According to Proposition A7.1 and Eq. (A12), $\forall x(t + t_0) \in W_t^1$, it holds that $x(t + t_0) \leq \sup W_{t_2-1}^1 < \frac{1}{4}pb + \frac{3}{4}gb$, thus the expression of $E_{x(t+t_0+1)}$ becomes (see Eq. (A7))

$$E_{x(t+t_0+1)} = \left[2x(t + t_0) - x(t + t_0 - 1), 2x(t + t_0) - x(t + t_0 - 1) + \frac{1}{4}C_2(gb - pb)\right]_I$$

(A13)

Therefore, we only need to prove that

$$\begin{cases} \sup_{x(t+t_0-1)\in W_{t-1}^1, x(t+t_0)\in W_t^1} [2x(t + t_0) - x(t + t_0 - 1)] \leq b_{t+1} \\ \inf_{x(t+t_0-1)\in W_{t-1}^1, x(t+t_0)\in W_t^1} \left[2x(t + t_0) - x(t + t_0 - 1) + \frac{1}{4}C_2(gb - pb)\right] \geq b_{t+1} + \delta_b \end{cases}$$

(A14)

Here,

$$\sup_{x(t+t_0-1)\in W_{t-1}^1, x(t+t_0)\in W_t^1} [2x(t + t_0) - x(t + t_0 - 1)] - b_{t+1}$$

$$\leq 2(b_t + \delta_b) - b_{t-1} - b_{t+1} = 0$$

and

$$\inf_{x(t+t_0-1)\in W_{t-1}^1, x(t+t_0)\in W_t^1} \left[2x(t + t_0) - x(t + t_0 - 1) + \frac{1}{4}C_2(gb - pb)\right] - (b_{t+1} + \delta_b)$$

$$\geq 2b_t - (b_{t-1} + \delta_b) + \frac{1}{4}C_2(gb - pb) - (b_{t+1} + \delta_b)$$

$$= \frac{1}{4}C_2(gb - pb) - 4\delta_b$$

$$\leq \frac{1}{4}d_0 + 4\frac{1}{16}d_0 = 0$$

□

By combining Propositions A7.1-A7.3, we successfully verify the feasibility of the state transition chain for agent's position $W_0^1 \to W_1^1 \to \cdots \to W_{t_2}^1$. According to Eq. (A12), it holds that
$$\inf W_{t_2}^1 = b_{t_2}$$
$$\geq \frac{1}{4}pb + \frac{3}{4}gb - \delta_b$$
$$\geq \frac{1}{4}pb + \frac{3}{4}gb - \frac{1}{16}(gb - pb)$$
$$> \frac{1}{2}(pb + gb)$$

This implies that once the position of an agent reaches $W_{t_2}^1$, the agent successfully reaches $U_m$. In summary, the feasibility of the state transition chain $\{\boldsymbol{\xi}(t_0)\} \to U_1^1 \to U_2^1 \to \cdots U_{t_2-1}^1 \to U_m$ for $\boldsymbol{\xi}(t_0) \in S_0$ satisfying $x(t_0 + 1) < \frac{1}{4}pb + \frac{3}{4}gb$ is successfully verified.

□

**A.8 Proof of Lemma 2(c)**

For $\boldsymbol{\xi}(t_0) \in S_0$ satisfying $x(t_0) \geq \frac{1}{4}pb + \frac{3}{4}gb$, according to Lemma 2(a) and Eq. (4.14), based on the state transition chain $\{\boldsymbol{\xi}(t_0)\} \to S_1^1 \to S_2^1 \to \cdots S_{t_1-1}^1 \to S_a$, the first hitting time from $\{\boldsymbol{\xi}(t_0)\}$ to $S_a$ has an upper bound $\left\lceil \frac{8(ub-lb)}{d_0} \right\rceil - 1$, and the corresponding transition probability has a lower bound (it is important to note that $\lambda = \frac{1}{4}$ in Lemma 2(a)):

$$p_{0a} \geq \prod_{t=1}^{t_1} h(|R_t^1|)$$
$$= h^{t_1}(\delta_a)$$
$$\geq h^{\left\lceil \frac{8(ub-lb)}{d_0} \right\rceil}\left(\min\left\{\frac{1-\lambda}{4}d_0, \frac{d_0^2}{128(ub-lb)}\right\}\right)$$
$$= h^{\left\lceil \frac{8(ub-lb)}{d_0} \right\rceil}\left(\min\left\{\frac{3}{16}d_0, \frac{d_0^2}{128(ub-lb)}\right\}\right)$$
$$= h^{\left\lceil \frac{8(ub-lb)}{d_0} \right\rceil}\left(\frac{d_0^2}{128(ub-lb)}\right)$$

For $\boldsymbol{\xi}(t_0) \in S_0$ satisfying $x(t_0) < \frac{1}{4}pb + \frac{3}{4}gb$, as described in Section 4.3, the state transition chain from $\{\boldsymbol{\xi}(t_0)\}$ to $S_a$ is constructed in two steps. In the first step, according to Lemma 2(b) and Eq. (4.14), the state transition chain from $\{\boldsymbol{\xi}(t_0)\}$ to $U_m = \{\boldsymbol{\xi}(t) | x(t) \geq \frac{1}{2}pb + \frac{1}{2}gb\}$ is constructed. The first hitting time from $\{\boldsymbol{\xi}(t_0)\}$ to $U_m$ has an upper bound $\left\lceil \frac{4(ub-lb)}{d_0} \right\rceil$, and the corresponding transition probability has a lower bound (denoted as $p_{S_0 \to U_m}$):

$$p_{S_0 \to U_m} \geq \prod_{t=1}^{t_2} h(|W_t^1|)$$

$$= h^{t_2}(\delta_b)$$
$$\geq h^{\left\lceil\frac{4(ub-lb)}{d_0}\right\rceil}\left(\frac{1}{16}d_0\right)$$

In the second step, given that $\inf W^1_{t_2} - \sup W^1_{t_2-1} > v(t_0) \geq \frac{1}{4}d_0$ (see the proof of Proposition A7.1) and $\inf W^1_{t_2} > \frac{1}{2}(pb+gb)$ (see the proof of Lemma 2(b)), it can be deduced that for the state vector $\boldsymbol{\xi}(t)$ within the state transition chain $\{\boldsymbol{\xi}(t_0)\} \to U^1_1 \to U^1_2 \to \cdots U^1_{t_2-1} \to U_m$ at iteration $t_2 + t_0$ (i.e. $\boldsymbol{\xi}(t_2+t_0) \in U_m$), both $x(t_2+t_0) \geq \frac{1}{2}pb + \frac{1}{2}gb$ and $v(t_2+t_0) \geq \frac{1}{4}d_0$ hold. Thus, according to Lemma 2(a), we can construct a transition chain from $\{\boldsymbol{\xi}(t_2+t_0)\}$ to $S_a$ starting at iteration $t_2+t_0$ (with $\lambda = \frac{1}{2}$ in Lemma 2(a)). According to Lemma 2(a) and Eq. (4.14), the first hitting time in this state transition chain has an upper bound $\left\lceil\frac{8(ub-lb)}{d_0}\right\rceil - 1$, and the corresponding transition probability has a lower bound (denoted as $p_{U_m \to S_a}$):

$$p_{U_m \to S_a} \geq h^{\left\lceil\frac{8(ub-lb)}{d_0}\right\rceil}\left(\min\left\{\frac{1-\lambda}{4}d_0, \frac{d_0^2}{128(ub-lb)}\right\}\right)$$
$$= h^{\left\lceil\frac{8(ub-lb)}{d_0}\right\rceil}\left(\min\left\{\frac{1}{8}d_0, \frac{d_0^2}{128(ub-lb)}\right\}\right)$$
$$= h^{\left\lceil\frac{8(ub-lb)}{d_0}\right\rceil}\left(\frac{d_0^2}{128(ub-lb)}\right)$$

By combining the two steps discussed above, for $\boldsymbol{\xi}(t_0) \in S_0$ satisfying $x(t_0) < \frac{1}{4}pb + \frac{3}{4}gb$, the first hitting time from $\{\boldsymbol{\xi}(t_0)\}$ to $S_a$ has an upper bound $\left\lceil\frac{4(ub-lb)}{d_0}\right\rceil + \left\lceil\frac{8(ub-lb)}{d_0}\right\rceil - 1 < 13\frac{ub-lb}{d_0}$, while the corresponding transition probability has a lower bound

$$p_{0a} \geq p_{S_0 \to U_m} p_{U_m \to S_a}$$
$$= h^{\left\lceil\frac{4(ub-lb)}{d_0}\right\rceil}\left(\frac{1}{16}d_0\right) h^{\left\lceil\frac{8(ub-lb)}{d_0}\right\rceil}\left(\frac{d_0^2}{128(ub-lb)}\right)$$
$$> h^{13\frac{ub-lb}{d_0}}\left(\frac{d_0^2}{128(ub-lb)}\right)$$

Finally, by combining the first hitting time and transition probability results for the two cases $x(t_0) \geq \frac{1}{4}pb + \frac{3}{4}gb$ and $x(t_0) < \frac{1}{4}pb + \frac{3}{4}gb$, we derive Eq. (4.18).

□

### A.9 Proof of Lemma 3

We initially set $t \geq 2$.
If $x(t-1) \geq \max\{pb, gb\}$, then it holds that:
$$v(t) = v(t-1) + C_1 r_1(pb - x(t-1)) + C_2 r_2(gb - x(t-1)) \leq v(t)$$
If $x(t-1) = lb$, then it holds that $v(t-1) \leq 0$ (otherwise $x(t-1) = B_{[lb,ub]_I}(x(t-$

$2) + v(t − 1)) \geq B_{[lb,ub]_I}(lb + v(t − 1)) > lb)$, therefore it holds that:
$$v(t) = v(t − 1) + C_1 r_1(pb − lb) + C_2 r_2(gb − lb) \leq (C_1 + C_2)(ub − lb)$$

If $lb < x(t − 1) < \max\{pb, gb\}$, then it holds that $v(t − 1) = x(t − 1) − x(t − 2)$, therefore it holds that:
$$v(t) = [x(t − 1) − x(t − 2)] + C_1 r_1(pb − x(t)) + C_2 r_2(gb − x(t))$$
$$\leq (C_1 + C_2 + 1)(ub − lb)$$

To sum up, for $t \geq 2$, it holds that
$$v(t) \leq \max\{v(1), (C_1 + C_2 + 1)(ub − lb)\}$$

Therefore, $\forall t \geq 0$, it holds that
$$v(t) \leq \max\{v(1), (C_1 + C_2 + 1)(ub − lb), v(0)\}$$
$$\leq \max\{v(0) + (C_1 + C_2)(ub − lb), (C_1 + C_2 + 1)(ub − lb)\}$$
$$= (C_1 + C_2)(ub − lb) + \max\{v(0), ub − lb\}$$
$$= (C_1 + C_2 + 1)(ub − lb)$$

where the last equality utilizes the condition $v(0) \leq ub − lb$, which holds naturally.

$\square$

**A.10 Proof of Lemma 4**

We initially prove the following Lemma A2:

*Lemma A2*

Given $\xi(t) = [ub, v(t)]$, and setting the distribution range of $v(t + 1)$ as $D_{x(t+1)}$, it holds that $D_{x(t+1)} \supset [v(t) − d_0, v(t)]$.

*Proof of Lemma A2*

According to Eq. (2.7), it holds that:
$$v(t + 1) = v(t) − C_1 r_1(ub − pb) − C_2 r_2(ub − gb) \quad (A15)$$

Setting $r_1 = r_2 = 0$ in Eq. (A15) yields $v(t + 1) = v(t)$. Conversely, setting $r_1 = 1, r_2 = 0$ in Eq. (A8) yields $v(t + 1) = v(t) − C_1(ub − pb) \leq v(t) − C_1(gb − pb) \leq v(t) − d_0$. Since $D_{x(t+1)}$ represents a closed interval (see Corollary 1(a)), it holds that $D_{x(t+1)} \supset [v(t) − d_0, v(t)]$.

$\square$

Returning to the proof of Lemma 4.

*Proof of Lemma 4(a)*

We further define $V_0^2 = \{v(t_0)\}$ and $V_{\left\lceil\frac{3v(t_0)}{2d_0}\right\rceil+1}^2 = [-\lambda d_0, -\mu d_0]$. According to Eq. (4.24), $V_t^2$ form a series of decreasing intervals (i.e. $\forall t = 0, \ldots, \left\lceil\frac{3v(t_0)}{2d_0}\right\rceil - 1, \sup V_{t+1}^2 < \inf V_t^2$), and it holds that $\inf V_{\left\lceil\frac{3v(t_0)}{2d_0}\right\rceil}^2 = 0$. It is important to note that $\forall t = 0, \ldots, \left\lceil\frac{3v(t_0)}{2d_0}\right\rceil - 1$,

$$\sup V_t^2 − \inf V_{t+1}^2$$
$$\leq v(t_0) − \left(t − \frac{1}{2}\right)\frac{v(t_0)}{\left\lceil\frac{3v(t_0)}{2d_0}\right\rceil} − \left[v(t_0) − (t + 1)\frac{v(t_0)}{\left\lceil\frac{3v(t_0)}{2d_0}\right\rceil}\right]$$
$$\leq \frac{3}{2}\frac{v(t_0)}{\left\lceil\frac{3v(t_0)}{2d_0}\right\rceil} \leq d_0$$

and

$$\sup V^2_{\left\lceil\frac{3v(t_0)}{2d_0}\right\rceil} - \inf V^2_{\left\lceil\frac{3v(t_0)}{2d_0}\right\rceil+1}$$

$$= v(t_0) - \left(\left\lceil\frac{3v(t_0)}{2d_0}\right\rceil - \frac{1}{2}\right)\frac{v(t_0)}{\left\lceil\frac{3v(t_0)}{2d_0}\right\rceil} - (-\lambda d_0)$$

$$\leq \frac{1}{2}\frac{v(t_0)}{\left\lceil\frac{3v(t_0)}{2d_0}\right\rceil} + \frac{1}{2}d_0$$

$$\leq \frac{1}{3}d_0 + \frac{1}{2}d_0 < d_0$$

Thus, according to Lemma A2, for any $t = 0, \ldots, \left\lceil\frac{3v(t_0)}{2d_0}\right\rceil$, and for any $v(t+t_0) \in V_t^2$, the distribution range of $v(t+t_0+1)$ covers $V_{t+1}^2$. Therefore, $\{\xi(t_0)\} \to S_1^2 \to S_2^2 \to \cdots S_{\left\lceil\frac{3v(t_0)}{2d_0}\right\rceil}^2 \to S_b(\lambda, \mu)$ is a valid state transition chain from $\{\xi(t_0)\}$ to $S_b(\lambda, \mu)$.

□

*Proof of Lemma 4(b)*

Consider the transition chain $\{\xi(t_0)\} \to S_1^2 \to S_2^2 \to \cdots S_{\left\lceil\frac{3v(t_0)}{2d_0}\right\rceil}^2 \to S_b(\lambda, \mu)$, the corresponding first hitting time has an upper bound $\left\lceil\frac{3v(t_0)}{2d_0}\right\rceil + 1 \leq 2\left\lceil\frac{v_u}{d_0}\right\rceil = t_{ab}$ (since $\frac{v_u}{d_0} \geq \frac{(C_1+C_2+1)(ub-lb)}{d_0} \geq 3$), and according to Eq. (4.22), the corresponding transition probability has a lower bound:

$$Pr \geq \prod_{t=0}^{\left\lceil\frac{3v(t_0)}{2d_0}\right\rceil} h(|V_{t+1}^2|)$$

$$= h\big((\lambda-\mu)d_0\big)h^{\left\lceil\frac{3v(t_0)}{2d_0}\right\rceil}\left(\frac{1}{2}\frac{v(t_0)}{\left\lceil\frac{3v(t_0)}{2d_0}\right\rceil}\right)$$

$$\geq h\big((\lambda-\mu)d_0\big)h^{t_{ab}}\left(\frac{1}{2}\frac{v(t_0)}{\frac{v(t_0)}{d_0}}\right)$$

$$= h\big((\lambda-\mu)d_0\big)h^{t_{ab}}\left(\frac{1}{4}d_0\right) = p_{ab}$$

□

**A.11 Proof of Lemma 5(a)**

We further define $R_{t_1}^3 = [\max\{a_{t_1} - \delta_a, l_g\}, a_{t_1}]$ and $R_0^3 = \{x(t_0)\}$. Obviously, it holds that $a_{t_1} = u_g$. we will sequentially introduce and prove the following two propositions.

*Proposition A11.1*

$R_t^3(t = 0, \ldots, t_1)$ represents a series of decreasing intervals. Specifically, it holds that $\forall t = 1, \ldots, t_1, \sup R_t^1 < \inf R_{t-1}^1$.

*Proof of Proposition A11.1*

For $t = 1$, it holds that $\inf R_{t-1}^1 - \sup R_t^1 = |v(t_0)|$; For any $t = 2, \ldots, t_1$, it holds that
$\inf R_{t-1}^1 - \sup R_t^1$
$= a_{t-1} - \delta_a - a_t$
$= |v(t_0)| + (2t-3)\delta_a$
$> |v(t_0)| > 0$

□

*Proposition A11.2*

For any $t = 1, \dots, t_1 - 1$, and for any $x(t + t_0 - 1) \in R_{t-1}^3, x(t + t_0) \in R_t^3$, the distribution range of $x(t + t_0 + 1)$ covers $R_{t+1}^3$.

*Proof of Proposition A11.3*

According to Lemma A1, we only need to prove that $\forall t = 1, \dots, t_1 - 1, \forall x(t + t_0 - 1) \in R_{t-1}^3, \forall x(t + t_0) \in R_t^3, E_{x(t+t_0+1)} \supset R_{t+1}^3$ holds, where $E_{x(t+t_0+1)}$ is defined in Eq. (A7). According to Proposition A11.1, $\forall t = 1, \dots, t_1 - 1, \forall x(t + t_0) \in R_t^3$, it holds that $x(t + t_0) > a_{t_1} = u_g \geq \frac{3}{4}pb + \frac{1}{4}gb$. Thus, the expression of $E_{x(t+t_0+1)}$ can be rewritten as follows (see Eq. (A7)):

$$E_{x(t+t_0+1)} = \left[2x(t+t_0) - x(t+t_0-1) - \frac{1}{4}C_1(gb - pb), 2x(t+t_0) - x(t+t_0-1)\right]_I$$

(A16)

Therefore, we only need to prove that

$$\begin{cases} \inf_{x(t+t_0-1) \in R_{t-1}^3, x(t+t_0) \in R_t^3} [2x(t+t_0) - x(t+t_0-1)] \geq a_{t+1} \\ \sup_{x(t+t_0-1) \in R_{t-1}^3, x(t+t_0) \in R_t^3} \left[2x(t+t_0) - x(t+t_0-1) - \frac{1}{4}C_1(gb - pb)\right] \leq a_{t+1} - \delta_a \end{cases}$$

(A17)

Here,

$$\inf_{x(t+t_0-1) \in R_{t-1}^3, x(t+t_0) \in R_t^3} [2x(t+t_0) - x(t+t_0-1)] - a_{t+1}$$

$$\geq 2(a_t - \delta_a) - a_{t-1} - a_{t+1} = 0$$

and

$$\sup_{x(t+t_0-1) \in R_{t-1}^3, x(t+t_0) \in R_t^3} \left[2x(t+t_0) - x(t+t_0-1) - \frac{1}{4}C_1(gb - pb)\right] - (a_{t+1} - \delta_a)$$

$$\leq 2a_t - (a_{t-1} - \delta_a) - \frac{1}{4}C_1(gb - pb)(gb - pb) - (a_{t+1} - \delta_a)$$

$$= -\frac{1}{4}C_1(gb - pb) + 4\delta_a$$

$$= 4\frac{ub - u_g - \left(\left|\frac{ub - u_g}{|v(t_0)|}\right| - 1\right)|v(t_0)|}{t_1(t_1 - 1)} - \frac{1}{4}d_0$$

$$\leq 8\frac{|v(t_0)|}{\left(\left|\frac{ub - u_g}{|v(t_0)|}\right| - 2\right)\left(\left|\frac{ub - u_g}{|v(t_0)|}\right| - 3\right)} - \frac{1}{4}d_0$$

$$\leq 8\frac{\frac{1}{20}d_0}{\left(\left|\frac{\frac{1}{3}d_0}{\frac{1}{20}d_0}\right| - 2\right)\left(\left|\frac{\frac{1}{3}d_0}{\frac{1}{20}d_0}\right| - 3\right)} - \frac{1}{4}d_0$$

$$= -\frac{13}{60}d_0 < 0$$

□

By combining Propositions A11.1-A11.2, the feasibility of the state transition chain $\{\xi(t_0)\} \to S_1^3 \to S_2^3 \to \cdots S_{t_1-1}^3 \to S_g$ is successfully validated.

□

## A.12 Proof of Lemma 5(b)

We further define $Y_{t_3}^3 = [\max\{c_{t_3}, l_g\}, c_{t_3} + \delta_c]$ and $Y_0^3 = \{x(t_0 + t_2)\}$. Obviously, it holds that $b_{t_2} = \frac{1}{4}pb + \frac{3}{4}gb$, $c_{t_3} + \delta_c = u_g$. For the initial state vector $\xi(t_0) \in S_b\left(\frac{1}{40}, \frac{1}{20}\right)$, similar to the proof of A.11, the feasibility of the state transition chain $\{\xi(t_0)\} \to U_1^3 \to U_2^3 \to \cdots U_{t_2-1}^3 \to U_{t_2}^3$ can be easily proven. Based on this, we will sequentially introduce and prove the following two propositions.

*Proposition A12.1*

$t_2 \geq 3$ and $\delta_b \leq \frac{1}{60}d_0$ hold.

*Proof of Proposition A12.1*

It holds that

$$t_2 = \left\lceil \frac{ub - \left(\frac{1}{4}pb + \frac{3}{4}gb\right)}{|v(t_0)|} \right\rceil - 2$$

$$\geq \left\lceil \frac{\frac{1}{4}(gb - pb)}{|v(t_0)|} \right\rceil - 2$$

$$\geq \left\lceil \frac{\frac{1}{4}d_0}{\frac{1}{20}d_0} \right\rceil - 2 = 3$$

and

$$\delta_b = \frac{ub - \left(\frac{1}{4}pb + \frac{3}{4}gb\right) - \left(\left\lceil \frac{ub - \left(\frac{1}{4}pb + \frac{3}{4}gb\right)}{|v(t_0)|} \right\rceil - 1\right)|v(t_0)|}{t_2(t_2 - 1)}$$

$$\leq \frac{2|v(t_0)|}{t_2(t_2 - 1)}$$

$$\leq \frac{2 \cdot \frac{1}{20}d_0}{6} = \frac{1}{60}d_0$$

□

*Proposition A12.2*

$|v(t_0 + t_2)| \leq \frac{8}{3}|v(t_0)|$ holds.

*Proof of Proposition A12.2*

It holds that

$|v(t_0 + t_2)| \leq \sup W_{t_2-1}^3 - \inf W_{t_2-1}^3$
$= b_{t_2-1} - (b_{t_2} - \delta_b)$
$= |v(t_0)| + (2t_2 - 1)\delta_b$
$\leq |v(t_0)| + (2t_2 - 1)\frac{2|v(t_0)|}{t_2(t_2-1)}$

$\leq |v(t_0)| + 5\frac{2|v(t_0)|}{6} = \frac{8}{3}|v(t_0)|$

□

*Proposition A12.3*

$t_3 \geq 5$ and $\delta_c \leq \frac{1}{75} d_0$ hold.

*Proof of Proposition A12.1*

It holds that

$$t_3 \geq \left\lceil \frac{\left(\frac{1}{4}pb + \frac{3}{4}gb - \delta_b\right) - \left(\frac{3}{4}pb + \frac{1}{4}gb\right)}{|v(t_0+t_2)|} \right\rceil + 1$$

$$= \left\lceil \frac{\frac{1}{2}(gb-pb) - \delta_b}{|v(t_0+t_2)|} \right\rceil + 1$$

$$\geq \left\lceil \frac{\frac{1}{2}d_0 - \frac{1}{60}d_0}{\frac{8\,1}{320}d_0} \right\rceil + 1 = 5$$

and

$$\delta_c = \frac{u_g - x(t_0+t_2) + \left(\left\lceil \frac{x(t_0+t_2) - u_g}{|v(t_0+t_2)|} \right\rceil + 1\right)|v(t_0+t_2)|}{t_3(t_3-1)}$$

$$\leq \frac{2|v(t_0+t_2)|}{t_3(t_3-1)}$$

$$\leq \frac{2 \frac{8\,1}{320}d_0}{20} = \frac{1}{75}d_0$$

□

*Proposition A12.4*

$Y_t^3 (t = 0, \ldots, t_3)$ represents a series of decreasing intervals. Specifically, it holds that $\forall t = 1, \ldots, t_3, \sup Y_t^3 < \inf Y_{t-1}^3$.

*Proof of Proposition A12.1*

For $t = 1$, it holds that $\inf Y_{t-1}^3 - \sup Y_t^3 = |v(t_0+t_2)| > 0$; For $t = 2, \ldots, t_3$, it holds that
$\inf R_{t-1}^1 - \sup R_t^1$
$= c_{t-1} - (c_t + \delta_c)$
$= |v(t_0 + t_2)| - (2t-1)\delta_c$
$\geq |v(t_0 + t_2)| - (2t_3 - 1) \frac{2|v(t_0+t_2)|}{t_3(t_3-1)}$
$\geq \frac{1}{10}|v(t_0 + t_2)| > 0$

□

*Proposition A12.5*

For any $t = 1, \ldots, t_3 - 1$, and for any $x(t + t_0 - 1) \in Y_{t-1}^3, x(t + t_0) \in Y_t^3$, the distribution range of $x(t + t_0 + 1)$ covers $Y_{t+1}^3$.

*Proof of Proposition A12.5*

According to Lemma A1, we only need to prove that $\forall t = 1, \ldots, t_3 - 1, \forall x(t + t_0 - 1) \in Y_{t-1}^3, \forall x(t + t_0) \in Y_t^3, E_{x(t+t_0+1)} \supset Y_{t+1}^3$ holds, where $E_{x(t+t_0+1)}$ is defined in Eq. (A7). According to Proposition A12.4, $\forall t = 1, \ldots, t_3 - 1, \forall x(t + t_0) \in Y_t^3$, it holds that $x(t + t_0) \leq \frac{1}{4}pb + \frac{3}{4}gb$. Thus, the expression of $E_{x(t+t_0+1)}$ can be rewritten as follows (see Eq. (A7)):

$$E_{x(t+t_0+1)} = \left[2x(t+t_0) - x(t+t_0-1), 2x(t+t_0) - x(t+t_0-1) + \frac{1}{4}C_2(gb-pb)\right]_I$$

(A18)

Therefore, we only need to prove that

$$\begin{cases} \sup_{x(t+t_0-1)\in Y_{t-1}^3, x(t+t_0)\in Y_t^3} [2x(t+t_0) - x(t+t_0-1)] \leq c_{t+1} \\ \inf_{x(t+t_0-1)\in Y_{t-1}^3, x(t+t_0)\in Y_t^3} \left[2x(t+t_0) - x(t+t_0-1) + \frac{1}{4}C_2(gb-pb)\right] \geq c_{t+1} + \delta_c \end{cases}$$

(A19)

Here,

$$\sup_{x(t+t_0-1)\in Y_{t-1}^3, x(t+t_0)\in Y_t^3} [2x(t+t_0) - x(t+t_0-1)] - c_{t+1}$$

$$\geq 2(c_t + \delta_c) - c_{t-1} - a_{t+1} = 0$$

and

$$\inf_{x(t+t_0-1)\in Y_{t-1}^3, x(t+t_0)\in Y_t^3} \left[2x(t+t_0) - x(t+t_0-1) + \frac{1}{4}C_2(gb-pb)\right] - (c_{t+1} + \delta_c)$$

$$\geq 2c_t - (c_{t-1} + \delta_a) + \frac{1}{4}C_2(gb-pb) - (c_{t+1} + \delta_c)$$

$$\geq \frac{1}{4}d_0 - 4\delta_c$$

$$\geq \frac{1}{4}d_0 - \frac{4}{75}d_0 > 0$$

□

By combining Propositions A12.1-A12.5, the feasibility of the state transition chain $\{\boldsymbol{\xi}(t_0 + t_2)\} \to X_1^3 \to X_2^3 \to \cdots X_{t_3-1}^3 \to S_g$ is successfully validated.

□

## A.13 Proof of Lemma 2(c)

For the case where $ub - u_g \leq \frac{1}{3}d_0$, it evidently holds that $S_b\left(\frac{ub-u_g}{d_0}, \min\left\{\frac{ub-l_g}{d_0}, \frac{1}{2}\right\}\right) \subset S_g$.

For the case where $ub - u_g > \frac{1}{3}d_0$ and $ug \geq \frac{3}{4}pb + \frac{1}{4}gb$, for $\boldsymbol{\xi}(t_0) \in S_b\left(\frac{1}{40}, \frac{1}{20}\right)$, according to Lemma 5(a) and Eq. (4.14), based on the state transition chain $\{\boldsymbol{\xi}(t_0)\} \to S_1^3 \to S_2^3 \to \cdots S_{t_1-1}^3 \to S_g$, the first hitting time from $\{\boldsymbol{\xi}(t_0)\}$ to $S_g$ has an upper bound $t_1 - 2 \leq \left\lceil \frac{ub-u_g}{\frac{1}{40}d_0} \right\rceil - 2 \leq \left\lceil 40\frac{ub-lb}{d_0} \right\rceil - 2$, and the corresponding transition probability has a lower bound:

$$p_{bg} \geq \prod_{t=1}^{t_1} h(|R_t^3|)$$
$$\geq h^{t_1}(\delta_a)h(|R_g|)$$
$$> h^{\left\lceil 40\frac{ub-lb}{d_0}\right\rceil - 2}\left(\frac{|v(t_0)|}{t_1(t_1-1)}\right)h(|R_g|)$$
$$> h^{\left\lceil 40\frac{ub-lb}{d_0}\right\rceil - 2}\left(\frac{\frac{1}{40}d_0}{\left(40\frac{ub-lb}{d_0}\right)^2}\right)h(|R_g|)$$
$$> h^{\left\lceil 40\frac{ub-lb}{d_0}\right\rceil - 2}\left(10^{-5}\frac{d_0^3}{(ub-lb)^2}\right)h(|R_g|)$$

For the case where $ub - u_g > \frac{1}{3}d_0$ and $ug < \frac{3}{4}pb + \frac{1}{4}gb$, for $\xi(t_0) \in S_b\left(\frac{1}{40}, \frac{1}{20}\right)$, and for the first step of the transition chain $\{\xi(t_0)\} \to U_1^3 \to U_2^3 \to \cdots U_{t_2-1}^3 \to U_{t_2}^3$, similarly, the first hitting time from $\{\xi(t_0)\}$ to $U_{t_2}^3$ has an upper bound $40\frac{ub-lb}{d_0}$, and the corresponding transition probability has a lower bound $h^{40\frac{ub-lb}{d_0}}\left(10^{-5}\frac{d_0^3}{(ub-lb)^2}\right)$. While for the second step of the transition chain $\{\xi(t_0+t_2)\} \to X_1^3 \to X_2^3 \to \cdots X_{t_3-1}^3 \to S_g$, the first hitting time from $\{\xi(t_0+t_2)\}$ to $S_g$ has an upper bound $\left\lceil\frac{x(t_0+t_2)-u_g}{|v(t_0+t_2)|}\right\rceil + 1 < \left\lceil\frac{ub-lb}{|v(t_0)|}\right\rceil + 1 \leq \left\lceil 40\frac{ub-lb}{d_0}\right\rceil + 1$, and the corresponding transition probability has a lower bound:

$Pr \geq \prod_{t=1}^{t_3} h(|Y_t^3|)$
$\geq h^{t_3}(\delta_c) h(|R_g|)$
$> h^{\left\lceil 40\frac{ub-lb}{d_0}\right\rceil+1}\left(\frac{|v(t_0+t_2)|}{t_3(t_3-1)}\right) h(|R_g|)$
$> h^{\left\lceil 40\frac{ub-lb}{d_0}\right\rceil+1}\left(\frac{\frac{1}{40}d_0}{\left\lceil 40\frac{ub-lb}{d_0}\right\rceil\left(\left\lceil 40\frac{ub-lb}{d_0}\right\rceil+1\right)}\right) h(|R_g|)$
$> h^{\left\lceil 40\frac{ub-lb}{d_0}\right\rceil+1}\left(10^{-5}\frac{d_0^3}{(ub-lb)^2}\right) h(|R_g|)$

By combining the two steps discussed above, the first hitting time from $\{\xi(t_0)\}$ to $S_g$ has an upper bound $\left\lceil 40\frac{ub-lb}{d_0}\right\rceil - 2 + \left\lceil 40\frac{ub-lb}{d_0}\right\rceil + 1 < 80\frac{ub-lb}{d_0}$, and the corresponding transition probability has a lower bound $h^{80\frac{ub-lb}{d_0}}\left(10^{-5}\frac{d_0^3}{(ub-lb)^2}\right) h(|R_g|)$.

Finally, by combining the first hitting time and transition probability results for the above two cases $ug \geq \frac{3}{4}pb + \frac{1}{4}gb$ and $ug < \frac{3}{4}pb + \frac{1}{4}gb$, we derive Eq. (4.18).

$\square$

**A.14 Proof of Proposition 2**

Obviously, for any initial state $\{\xi(t_0)\}$, to determine the values of $t_{e0}$ and $p_{e0}$, we only need to analyze the state transition chain $\{\xi(t_0)\} \to S_0 \to \cdots \to S_a \to \cdots \to S_b\left(\frac{1}{40}, \frac{1}{20}\right) \to \cdots \to S_g$. According to Eq. (4.11) and Lemmas 2, 3, and 5, the corresponding first hitting time has an upper bound:

$t \leq 1 + t_{0a} + t_{ab} + t_{bg}$
$= 1 + 13\left\lceil\frac{ub-lb}{d_0}\right\rceil + 2\left\lceil\frac{v_u}{d_0}\right\rceil + \left\lceil 80\frac{ub-lb}{d_0}\right\rceil$
$\leq 5 + 13\frac{ub-lb}{d_0} + 2\frac{(C_1+C_2+1)(ub-lb)}{d_0} + 80\frac{ub-lb}{d_0}$
$\leq t_{e0}$

The corresponding transition probability has a lower bound:

$pr \geq h\left(\frac{1}{4}d_0\right) p_{0a} p_{ab} p_{bg}$

$$\geq h\left(\frac{1}{4}d_0\right) h^{t_{oa}}\left(\frac{d_0^2}{128(ub-lb)}\right) h\left(\frac{1}{40}d_0\right) h^{t_{ab}}\left(\frac{1}{4}d_0\right) h^{t_{bg}}\left(10^{-5}\frac{d_0^3}{(ub-lb)^2}\right) h(|R_g|)$$

$$\geq h^{2+t_{oa}+t_{ab}+t_{bg}}\left(\min\left\{\frac{1}{4}d_0, \frac{1}{40}d_0, \frac{d_0^2}{128(ub-lb)}, 10^{-5}\frac{d_0^3}{(ub-lb)^2}\right\}\right) h(|R_g|)$$

$$\geq h^{t_{eo}}\left(10^{-5}\frac{d_0^3}{(ub-lb)^2}\right) h(|R_g|) = p_{e0}$$

□

The necessary and sufficient conditions for PSO to ensure that the probability of escaping from LO within a finite number of iterations is equal to 1 are derived. On this basis, two behaviors of agents, namely inertial motion and oscillation, are proposed, and the state transition chain leading towards the global optimum is constructed.